\documentclass[priprint,12pt,A4paper,3p]{elsarticle}

\usepackage{amssymb}
\usepackage{amsmath}

\usepackage{lipsum,bm}
\usepackage{amsfonts,wasysym}
\usepackage{graphicx}
\usepackage{epstopdf}
\usepackage{algorithmic}
\usepackage{amsopn}
\usepackage{amssymb}
\usepackage{mathtools}
\usepackage{verbatim}
\usepackage{enumitem}
\usepackage{mathrsfs}
\usepackage{multirow}
\usepackage{subfig}
\usepackage{caption}   
\usepackage{float} 
\usepackage{booktabs}
\usepackage{diagbox}
\usepackage{hyperref}
\ifpdf
\DeclareGraphicsExtensions{.eps,.pdf,.png,.jpg}
\else
\DeclareGraphicsExtensions{.eps}
\fi

\usepackage{amsmath, amsthm}
\newtheorem{Corollary}{Corollary}
\newtheorem{Assumption}{Assumption}
\newtheorem{theorem}{Theorem}
\newtheorem{Proposition}{Proposition}
\newtheorem{lemma}{Lemma}

\theoremstyle{definition}
\newtheorem{Definition}{Definition}

\newcommand{\mS}{\mathcal{S}}
\newcommand{\norm}[1]{\left\Vert#1\right\Vert}
\newcommand{\abs}[1]{\left\vert#1\right\vert}
\newcommand{\set}[1]{\left\{#1\right\}}
\newcommand{\seq}[1]{\left<#1\right>}

\newcommand{\E}{\mathbb{E}}

\newcommand{\R}{\mathbb{R}}
\newcommand{\N}{\mathbb{N}}
\newcommand{\PP}[1]{\mathbb{P}\left[#1\right]}
\newcommand{\diag}[1]{\mathrm{diag}(#1)}

\newcommand{\x}{\boldsymbol{x}}
\newcommand{\X}{\mathcal{X}}
\newcommand{\y}{\boldsymbol{y}}
\newcommand{\uu}{\boldsymbol{u}}
\newcommand{\vv}{\boldsymbol{v}}

\newcommand{\dd}{\mathrm{d}}
\newcommand{\cU}{\mathcal{U}}
\newcommand{\cL}{\mathcal{L}}
\newcommand{\cE}{\mathcal{E}}
\newcommand{\hcL}{\widehat{\mathcal{L}}}
\newcommand{\cZ}{\mathcal{Z}}
\newcommand{\cV}{\mathcal{V}}
\newcommand{\hcU}{\widehat{\mathcal{U}}}

\newcommand{\hZ}{\widehat{Z}}
\newcommand{\bhZ}{\overline{\widehat{Z}}}

\newcommand{\eexp}[1]{\exp\left\{#1\right\}}

\newcommand{\marker}[1]{{\color{red}#1}}

\begin{document}

\begin{frontmatter}



\title{Generalization Error Analysis of Deep Backward Dynamic Programming for Solving Nonlinear PDEs} 


\author[label1]{Du Ouyang\corref{cor1}}
\ead{oyd21@mails.tsinghua.edu.cn}
\author[label1]{Jichang Xiao}\ead{xiaojc19@mails.tsinghua.edu.cn}
\author[label1]{Xiaoqun Wang}\ead{wangxiaoqun@mail.tsinghua.edu.cn}
\affiliation[label1]{organization={Department of Mathematical Sciences, Tsinghua University},
            city={Beijing},
            postcode={100084}, 
            country={People's Republic of China}}
\cortext[cor1]{Corresponding author.}
\begin{abstract}
We explore the application of the quasi-Monte Carlo (QMC) method in deep backward dynamic programming (DBDP)~\cite{hure2020} for numerically solving high-dimensional nonlinear partial differential equations (PDEs). Our study focuses on examining the generalization error as a component of the total error in the DBDP framework, discovering that the rate of convergence for the generalization error is influenced by the choice of sampling methods. Specifically, for a given batch size $m$, the generalization error under QMC methods exhibits a convergence rate of $O(m^{-1+\varepsilon})$, where $\varepsilon$ can be made arbitrarily small. This rate is notably more favorable than that of the traditional Monte Carlo (MC) methods, which is $O(m^{-1/2+\varepsilon})$. Our theoretical analysis shows that the generalization error under QMC methods achieves a higher order of convergence than their MC counterparts. Numerical experiments demonstrate that QMC indeed surpasses MC in delivering solutions that are both more precise and stable.
\end{abstract}


\end{frontmatter}



\section{Introduction}
	Nonlinear parabolic partial differential equations (PDEs) represent a general form of linear Kolmogorov PDEs~\cite{beck2021solving,berner2020numerically,de2022error,glau2022deep,richter2022robust}, which find wide applications in  financial mathematics~\cite{karoui97}. The nonlinear PDE has the following form 
	\begin{equation} \label{eq:PDE}
		\left\{
		\begin{aligned}
			\frac{\partial u}{\partial t} + \mathscr{A} u + f(t,x,u,\Sigma^{\top} D_x u) & = 0 , \;\;\;\;\;\; \text{ for  } (t,x)\in [ 0,T)\times\R^d, \\
			u(T,x) &=g(x), \;\;\;\;\;  \mbox{ for } x \in \R^d,
		\end{aligned}
		\right.
	\end{equation}
	where $ f :[0,T]\times \R^d \times \R \times \R^{d} \rightarrow \R^d$ and $ g :\R^d \rightarrow \R $ is called the terminal function. The second-order generator $\mathscr{A}$ has the form
	\begin{equation}\label{eq:derivative_operator}
		\mathscr{A} := \frac{1}{2} \sum_{i,j =1}^d (\Sigma \Sigma^{\top})_{ij}\frac{\partial^2}{\partial x_i\partial x_j} + \sum_{i = 1}^d M_i \frac{\partial}{\partial x_i},
	\end{equation}
	where $(M_1(t,x),\dots,M_d(t,x))\coloneqq M(t,x)$ is a function defined on $[0,T] \times \R^d$ with values in $\R^d$, $\Sigma(t,x)$ is a function defined on $[0,T] \times \R^d$ with values in  $\mathbb{M}^d$, the set of  $d \times d$ matrices.
	
	The objective of this paper is to obtain the solution of high-dimensional PDE~\eqref{eq:PDE} at $t=0$ over a region $[a,b]^d$, i.e., we aim to obtain $u(0,x)$ for $x\in [a,b]^d$. Thanks to the nonlinear Feynman-Kac formula obtained by Pardoux and Peng~\cite{pardoux92}, the solution of the PDE~\eqref{eq:PDE} is connected to the backward stochastic differential equations (BSDEs). The BSDE was first proposed by Pardoux and Peng~\cite{pardoux90}, and has the following form
	\begin{equation}\label{eq:bsde}
		Y_t = g(\X_T) + \int_t^T f(s,\X_s,Y_s,Z_s) \dd s - \int_t^T Z_s^{\top} \dd W_s\ , \ \ \ 0\le t\le T,
	\end{equation}
	where $f$ and $g$ are identical to those in PDE~\eqref{eq:PDE}. The $\{\X_t\}_{0\le t\le T}$ satisfies the following forward stochastic differential equation (SDE)
	\begin{equation}\label{eq:SDE}
		\X_t = \eta + \int_0^t M(s,\X_s) \dd s+  \int_0^t  \Sigma(s,\X_s)  \dd W_s, \;\;\; 0 \le t \le T ,
	\end{equation}
	where  $W$ is a $d$-dimensional Brownian motion and $\eta $ is uniformly distributed on $\left[ a,b\right]^d$. Consider the probability space $(\Omega,\mathcal{F},P)$ equipped with a filtration 
	$ \{\mathcal{F}_t\}_{0\le t\le T}$ generated by the Brownian motion $ W$ and the initial random variable $\eta$ satisfying the usual conditions. $M$ and $\Sigma $ are idential to those in~\eqref{eq:PDE}.

	Under appropriate conditions (see~\cite{pardoux92} and~\cite{karoui97}), the solutions of~\eqref{eq:bsde} and~\eqref{eq:PDE} are connected in the form of
	\begin{equation}\label{eq:Y}
		Y_t = u(t,\X_t),\ \ 0\le t\le T,
	\end{equation}
	and if there are more smooth conditions, then
	\begin{equation}\label{eq:Z}
		Z_t = \Sigma^{\top}(t,\X_t)\partial_xu(t,\X_t),\ \ 0\le t \le T.
	\end{equation}
	The greatest advantage of the Feynman-Kac formula is its ability to transform the solution of PDEs into the solution of BSDEs. This is beneficial because general methods for solving PDEs, such as finite element methods, depend on the partition of the state space and are subject to the “curse of dimensionality”. The schemes in~\cite{bouchard04,gobet05,zhang04}, using the time partition and regression, overcome the issue of dimensionality. These schemes can be applied up to dimension 6 or 7. However, the efficacy of the regression methods is constrained by the choice of the basis functions (e.g., polynomial functions). In cases of higher dimensions (e.g., 50), the sheer volume of necessary basis functions renders these approaches highly inefficient.

	To address this issue, the pioneering deep BSDE method~\cite{e2017,han2018solving} applies deep neural networks to efficiently solve high-dimensional BSDEs (nonlinear PDEs). This approach  utilizes the required initial values as parameters and employs feedforward neural networks to approximate each time step to solve a global optimization problem.  In contrast, the Deep Backward Dynamic Programming (DBDP)~\cite{hure2020} simultaneously estimates the solution and its gradient through the minimization of sequential loss functions via backward induction, showing superior stability and accuracy compared to~\cite{e2017}, as evidenced by several test cases. Based on~\cite{hure2020}, Germain et al.~\cite{germain2022approximation} proposed a multistep version of DBDP. The deep splitting (DS) method~\cite{beck2021Deep} estimates the PDE solution via backward explicit local optimization problems, leveraging a neural network regression method to calculate conditional expectations. 

    In the aforementioned methods, the loss function is typically represented by an expectation of the states, which cannot be precisely calculated. To address this issue, empirical risk minimization (ERM) principle aims to apply the Monte Carlo (MC) method to sample states and using the sample mean (also called empirical risk) to estimate the loss function. Then the problem becomes finding a parameter that minimizes the empirical risk. Due to the ERM, the total error of the algorithm will be influenced by the generalization error (see Section~\ref{sec:total_error}). The generalization error (also called stochastic error in \cite{jiao2024error,jiao2023deep}), which correlates with the sampling method and sample size, signifies the capacity of ERM. As a component of the total error, it affects the performance of algorithm during the learning phase. Therefore, studying the convergence order of the generalization error, which is not involved in~\cite{beck2021Deep,germain2022approximation,hure2020}, is also important. Other research concerning the generalization error can be found in \cite{beck2022full,berner2020analysis,de2022error,jentzen2023overall,jiao2024error,jiao2023deep,mishra2023estimates}.

    The quasi-Monte Carlo (QMC) method is an efficient numerical method widely used in financial engineering (see~\cite{glasserman2004,hull2003}). For $ m$ quadrature points, the QMC method has a convergence rate $ O(m^{-1+\varepsilon})$ with an arbitrarily small $ \varepsilon > 0$, which is asymptotically faster than the MC rate $O(m^{-1/2})$. The applications of QMC methods in machine learning have also made some progress. Dick and Feischl~\cite{dick2021quasi} applied QMC methods to reduce the computational effort required for data compression in machine learning. Liu and Owen~\cite{liu2021quasi} demonstrated that randomized quasi-Monte Carlo (RQMC) methods combined with stochastic quasi-Newton methods can significantly accelerate the optimization. Longo et al.~\cite{longo2021higher} and Mishra et al.~\cite{mishra2021enhancing} utilized QMC methods to generate training data for training neural networks. Their results indicate that the QMC-based deep learning algorithm is more efficient than the MC-based one in approximating data-to-observables maps. Thus it is natural to ask whether the QMC methods can achieve higher efficiency compared to MC in solving nonlinear PDEs. Xiao et al~\cite{xiao2023analysis} studied the efficiency of RQMC methods for solving linear Kolmogorov PDEs by using deep learning. However, to the best of our knowledge, there is no related research for the nonlinear PDEs cases.
    
    To this end, we decompose the total error of the DBDP method for solving nonlinear PDEs into three parts: the first part is the scheme error, which relates to the time partition; the second part is the approximation error, related to the structure of the network; and the third part is the generalization error (ignored by~\cite{hure2020}), which is related to the sampling method and batch size. We study the convergence rate of the generalization error under MC and QMC frameworks. Our theoretical findings demonstrate that under QMC, the convergence rate of generalization error is $ O(m^{-1+\varepsilon})$ with $ \varepsilon>0$ being arbitrarily small, which is better than the MC case $ O(m^{-1/2+\varepsilon})$. The efficacy of the QMC method is validated through a series of experiments. The numerical results reveal that the QMC method yields training outcomes with significantly reduced relative error in comparison to the MC method when evaluated against the analytical solution. Moreover, the results obtained by the QMC method exhibit enhanced stability, characterized by a notably lower variance.

    The structure of this paper is as follows. Section~\ref{sec:pre} introduces the neural network and the DBDP scheme. Section~\ref{sec:total_error} decomposes the total error of the DBDP scheme and derives the generalization error. Section~\ref{sec:main} studies the convergence rates of generalization error under the MC and QMC frameworks. Section~\ref{sec:numerical} presents the numerical result of several examples. Section~\ref{sec:conclusion} concludes this paper.

	\section{Solving nonlinear PDE by using deep learning}\label{sec:pre} 
 Consider the time grid $\pi = \{(t_0,\dots,t_N) : 0= t_0<t_1<\dots<t_N=T\}$ with $|\pi| = \max_{0\le i\le N-1} |t_{i+1} - t_i|$. Recall the nonlinear Feynman-Kac formula~\eqref{eq:Y} and~\eqref{eq:Z}. The approximation of the solution $Y_{t_i}$ and $Z_{t_i}$ boils down to two parts. First, we use $X_{t_i}$ (this can be obtained by Euler-Maruyama or Milstein methods, see~\cite{kloeden1992,milstein06}) to estimate $\X_{t_i}$. Second, we approximate the function $u(t,\cdot)$ by neural network $\cU$.

	\subsection{Neural network}
	
	\begin{Definition}[Feedforward Neutral Network]\label{defn:fnn}
		Let $d,d_1 \in \mathbb{N}$, a feedforward neural network is a function of the form
		\begin{equation}\label{eq:mlp_network}
			x\in \R^d \mapsto A_k\circ\rho \circ A_{k-1} \circ \cdots \circ \rho \circ A_1(x) \in \R^{d_1},
		\end{equation}
		where $\rho : \R \rightarrow \R$ called activation function is a nonlinear continuous function and applied component-wise on the vectors, i.e., $\rho(x_1,\dots,x_m) = (\rho(x_1),\dots,\rho(x_m))$, and every $A_i: \R^{l_i}\rightarrow \R^{l_{i+1}}, i = 1,\dots,L$, is an affine transformation given by
		\begin{equation}
			A_i(x) = Q_ix + b_i.
		\end{equation}

  In the above equation, $Q_i$ and $ b_i$ are parameters called weight and bias, respectively. Denote $\theta = (Q_i,b_i)_{i = 1}^k$ be the collections of weight and bias, we can represent the functions in~\eqref{eq:mlp_network} by $\mathcal{U}(x;\theta)$. 
	\end{Definition}

In this paper, we choose the activation function $\rho$ to be the hyperbolic tangent function, i.e., 
		\begin{equation}\label{eq:tanh}
			\rho(x) = \frac{e^{x}-e^{-x}}{e^{x}+e^{-x}}.
		\end{equation}
 
		Let $\mathcal{S}=(l_1,\cdots,l_{k+1})$ be the structure of the network with $l_1 = d$ and $l_{k+1} = d_1$. 
		Denote $D(\mathcal{S})=k$ be the depth of the network, $\norm{\mathcal{S}}_\infty=\max_{1\leq i \leq k+1} l_i$ be the width of the network and $|\mathcal{S}| = \sum_{i = 1}^k (l_i\times l_{i+1} + l_{i+1})$ be the total number of the parameters of the network. We consider the parameter space with bound $ R$ and denote $ \Theta := \{\theta\mid \|\theta\|_{\infty} \le R\}$. Denote the class of neural networks \eqref{eq:mlp_network} with the same structure $\mathcal{S}$  by $\mathcal{N}_{\mathcal{S},R}=\{\mathcal{U}(x;\theta)\mid \theta \in \Theta\}$.

	The next lemma states some basic properties of the neural networks in Definition \ref{defn:fnn}. Similar results can be found in \cite{berner2020analysis,jiao2024error}.
	\begin{lemma}\label{lem:property_FNN}
		Under Definition \ref{defn:fnn}, suppose that $l_{k+1}=1$, for every $\theta_{1},\theta_2 \in \Theta $ and every $x\in \mathbb{R}^{d}$, we have
		\begin{equation}
			\abs{\mathcal{U}(x;\theta_1)-\mathcal{U}(x;\theta_2)} \leq C_{\mS,R}\max\{\norm{x}_\infty,1 \}\norm{\theta_1-\theta_2}_\infty,\nonumber 
		\end{equation}
		and
		\begin{equation}
			\abs{\mathcal{U}(x;\theta_1)} \leq B_{\mS,R},\nonumber
		\end{equation}
		where $C_{\mS,R}=|\mS|R^{D(\mS)-1}\norm{\mS}_\infty^{D(\mS)-1}$ and $B_{\mS,R}=2R\norm{\mS}_\infty$.
	\end{lemma}
	\begin{proof}
	    This lemma is a natural extension of \cite[Lemma 5.11]{jiao2024error}.
	\end{proof}
	
	\subsection{DBDP schemes}
	We consider the estimate of $Y_{t_i}$ as follows
	\begin{align}
		Y_{t_i} \approx \mathcal{U}_i(X_{t_{i}}),
	\end{align}
	where $\mathcal{U}_i$ is a neural network having the form~\eqref{eq:mlp_network} and $X_{t_{i}}$ is the approximation of $\X_{t_i}$. If SDE~\eqref{eq:SDE} does not have an analytical solution, then $X_{t_i}$ can be obtained by some numerical schemes, such as the Euler-Maruyama method and the Milstein method (see~\cite{kloeden1992,milstein06}).
	
	The deep backward dynamic programming (DBDP) scheme (the DBDP1 scheme in~\cite{hure2020}) is training $\mathcal{U}_i$ backward in time and  approximating $Z_{t_i}$ by another neural network $\mathcal{Z}_i(X_{t_i})$. Assume that
	\begin{equation}\label{eq:F}
		F(t, x, y, z, h, \Delta):=y-f(t, x, y, z) h + z^{\top} \Delta. 
	\end{equation}
	For functions $\hcU$, $\cU$ and $\cZ$, define
	\begin{equation}\label{eq:def_L_i} 
		\begin{aligned}
			L_i (\hcU,\cU,\cZ):= \E\big{|} \widehat \cU (X_{t_{i+1}}) - F(t_i,X_{t_i},\cU(X_{t_i}),\cZ(X_{t_i}),\Delta t_i,\Delta W_{t_i})\big{|}^2.
		\end{aligned}
	\end{equation}
	Assume $\mS_1 = (d,l_2,\dots,l_k,1)$ and $\mS_2 = (d,l_2,\dots,l_k,d)$. DBDP has the following form,\label{em:DBDP1}
	\begin{itemize}
		\item Initialize $\widehat \cU_N =g$.
		\item For $i = N-1,\dots,0$, given $\widehat \cU_{i+1}$, use $(\cU_i(\cdot,\theta),\cZ_i(\cdot,\theta)) \in \mathcal{N}_{\mS_1,R}\times \mathcal{N}_{\mS_2,R}$ for the approximation of $(u(t_i,\cdot),\Sigma^{\intercal}(t_i,\cdot)\partial_xu(t_i,\cdot))$, and find (for example, by SGD) the optimal $\theta$ to minimize the loss function
	\end{itemize}
	\begin{equation}\label{agm:DBDP1}
		\mathcal{L}_i(\theta) = L_i(\hcU_{i+1},\cU(\cdot;\theta),\cZ(\cdot;\theta)).
	\end{equation}
	\begin{itemize}
		\item Find the optimal $ \theta_i^* \in \arg\min_{\theta} \mathcal{L}_i(\theta)$ and update $\widehat \cU_i = \cU_i(\cdot;\theta_i^*)$.
	\end{itemize}

	\section{Total error of DBDP}\label{sec:total_error}
	In this section, we study the total error of the DBDP scheme. Recall that $Y = (Y_{t_0},\dots,Y_{t_{N-1}})$ and $Z = (Z_{t_0},\dots,Z_{t_{N-1}})$ are the solution of BSDE~\eqref{eq:bsde}. Let $\cU = (\cU_0,\dots,\cU_{N-1})$ and $\cZ = (\cZ_0,\dots,\cZ_{N-1})$ be two lists of neural networks. Define the error between $(\cU,\cZ)$ and $(Y,Z)$ as
	\begin{equation}\label{eq:total_error}
		\begin{aligned}
			\cE[(\cU,\cZ),(Y,Z)] = &\max_{i = 0,\dots,N-1} \E\left| Y_{t_i} - \cU_i(X_{t_i})\right|^2 \\
			&+ \E\left[ \sum_{i =0}^{N-1} \int_{t_i}^{t_{i+1}}\left| Z_t - \cZ_i(X_{t_i})\right|^2\dd t\right].
		\end{aligned}
	\end{equation}
	We make the following usual assumptions on the BSDE~\eqref{eq:bsde}. 
	\begin{Assumption}\label{assump:1}
		Assume that $f,g,M,\Sigma$ satisfy the lipschitz condition, i.e., there are constants $C_f, C_g, C>0$, such that for all $(t_1,x_1,y_1,z_1)$ and $(t_2,x_2,y_2,z_2) \in [0,T]\times\R^d\times\R\times\R^d$, 
		\begin{equation}\notag
			\left\{\begin{aligned}
				|f(t_1,x_1,y_1,z_1) - &f(t_2,x_2,y_2,z_2)| \\
				\le& C_f \left( |t_1-t_2|^{1/2} + \|x_1 - x_2\|_2 + |y_1 - y_2| + \|z_1 - z_2\|_2\right),\\
				|g(x_1) - g(x_2)| \le& C_g\|x_1-x_2\|_2,
			\end{aligned}\right.
		\end{equation}
		and
		\begin{equation}\notag
			\begin{aligned}
				|M(t_1,x_1) - M(t_2,x_2)| +  |\Sigma(t_1,x_1) - \Sigma(t_2,x_2)|
				\le C \left( |t_1-t_2|^{1/2} + \|x_1 - x_2\|_2 \right).
			\end{aligned}
		\end{equation}
	\end{Assumption}
	\begin{Assumption}\label{assump:2}
		Assume that the approximation $\{X_{t_i}\}$ of $\{\X_{t_i}\}$ satisfies
		\begin{equation}\notag
			\max_{i = 0,\dots,N-1}\E\left[ |X_{t_{i+1}} - \X_{t_{i+1}}|^2\right] = O(|\pi|),
		\end{equation}
            where $ |\pi| = \max_{i=0,\dots,N-1}|t_{i+1}-t_i|$.
	\end{Assumption}
 
	The following theorem can be proved by using the method in Theorem 4.1 of Hure et al.~\cite{hure2020}.
	\begin{theorem}\label{thm:total_error}
		Under Assumptions~\ref{assump:1} and~\ref{assump:2}, for functions $\cU = (\cU_0,\dots,\cU_{N-1})$ and $\cZ = (\cZ_0,\dots,\cZ_{N-1})$, there is a constant $C$, depending only on the coefficients of BSDE~\eqref{eq:bsde}, such that the error~\eqref{eq:total_error} has the following upper bound
		\begin{equation}\label{eq:total_error_final}
			\begin{aligned}
				\cE\left[ \left( \cU,\cZ\right),\left( Y,Z\right)\right] \le C\bigg(  |\pi| 
				+ N\sum_{i =0}^{N-1}\Big( L_i\left( \cU_{i+1},\cU_i,\cZ_i\right) - L_i\left(\cU_{i+1},\cU_i^*,\cZ_i^*\right)\Big)\bigg).
			\end{aligned}
		\end{equation}
		where $\cU_{N} = g$ and $\left( \cU_i^*,\cZ_i^*\right) = \arg\min_{(\cU_i,\cZ_i)} L_i\left( \cU_{i+1},\cU_i,\cZ_i\right)$.
	\end{theorem}
	\begin{proof}
		See~\ref{appendix:total_error}.
	\end{proof}
	
	If we fix a time partition $\pi$, then the total error is determined by the second term on the right hand side of~\eqref{eq:total_error_final}. Our goal is to find $\cU$ and $\cZ$ such that each term $L_i(\cU_{i+1},\cU_i,\cZ_i) - L_i(\cU_{i+1},\cU_i^*,\cZ_i^*)$ is as small as possible, enabling $(\cU,\cZ)$ to approximate $(Y,Z)$ with the minimal error. Suppose that $\cU_i$ and $\cZ_i$ are neural networks with parameter $\theta_i \in \Theta$. We have the following corollary, which can be proved by the method in the~\ref{appendix:total_error}.
	
	\begin{Corollary}\label{cor:error_decom}
		Under the conditions of Theorem~\ref{thm:total_error}, suppose that 
 $\cU_i = \cU_i(\cdot;\theta_i)$ and $\cZ_i = \cZ_i(\cdot;\theta_i)$, then
		\begin{equation}\label{eq:L_decomposition}
			\begin{aligned}
				L_i(&\cU_{i+1},\cU_i,\cZ_i) - L_i(\cU_{i+1},\cU_i^*,\cZ_i^*)\\
				\le& \underbrace{\inf_{\theta \in \Theta} \E\left[ \left| \cU_i^*(X_{t_i}) - \cU_{i}(X_{t_i};\theta)\right|^2\right] +\inf_{\theta \in \Theta} \E\left[ \left| \cZ_i^*(X_{t_i}) - \cZ_{i}(X_{t_i};\theta)\right|^2\right]}_{\text{Approximation Error}} \\
				&+  \mathcal{L}_i(\theta_i) -  \mathcal{L}_i(\theta_i^*),
			\end{aligned}
		\end{equation}
		where $\cU_i^*$ and $ \cZ_i^*$ are two functions in Theorem~\ref{thm:total_error} and $$ \mathcal{L}_i(\theta_i) := L_i(\cU_{i+1}(\cdot;\theta_{i+1}),\cU_i(\cdot;\theta_i),\cZ_i(\cdot;\theta_i)), \quad \theta^*_i \in \arg\min_{\theta_i\in\Theta} \mathcal{L}_i(\theta_i).$$
	\end{Corollary}
 
	The first two terms on the right hand side of~\eqref{eq:L_decomposition} have the same form as in Theorem 4.1 of Hure et al.~\cite{hure2020}. We refer to the sum of these two terms as the approximation error since it quantifies how accurately the functions can be approximated by the neural networks. The approximation error can be arbitrarily small if the network and $\Theta$ are large enough, see \cite{hornik1991approximation,hornik1993some,guhring2020error,siegel2020approximation,elbrachter2021deep,grohs2023proof,jentzen2021proof,hutzenthaler2020proof,hutzenthaler2020overcoming,darbon2016algorithms,hutzenthaler2023overcoming}.   The last term quantifies how well we can optimize the loss function to obtain the optimal parameters $\theta_i^{*}$. This term vanishes if we choose $\theta_i$ to be $\theta^*_i$ for $i = N-1,\dots,0$ (this is what Hure et al.~\cite{hure2020} do). 
 
    However, we cannot obtain $\theta^*_i$ by minimizing  $\mathcal{L}_i(\theta_i)$ directly, since the expectation in~\eqref{eq:def_L_i} is usually intractable. There are various methods to optimize the loss functional formulated as an expectation, including the Stochastic Gradient Descent~(SDG) method \cite{robbins1951stochastic,kingma2014adam,ruder2016overview} and 
 the Empirical risk minimization~(ERM) method \cite{vapnik1991principles,shalev2014understanding,mohri2018foundations}. In those methods, the true risk $\mathcal{L}_i$ is usually approximated by the empirical risk $\mathcal{L}_i^{(m)}$, which is exactly the sample mean over $m$ quadrature points.   Typically, we use the MC method to simulate the points, but the QMC method (see subsection~\ref{sec:qmc}) is an alternative quadrature method which achieves smaller errors when approximating expectations. From this prospective, the optimal tractable parameter is the minimizer of the empirical risk $\mathcal{L}_i^{(m)}$, which can be obtained by ERM.  Denote the minimizer of the empirical risk computed by MC and QMC methods by $\theta_i^{(m)}$ and $\widehat\theta_i^{(m)}$, respectively. To study the impact of the sampling methods, we investigate the error between $\theta_i^{(m)}$ (or $\widehat\theta_i^{(m)}$) and the global optimal parameter $\theta^*_i$, which is defined as the following (under MC, for illustration),
	\begin{equation}\label{eq:generalization_error}
		\mathcal{L}_i(\theta_i^{(m)}) - \mathcal{L}_i(\theta_i^*).
	\end{equation}
    We follow \cite{beck2022full} to refer to \eqref{eq:generalization_error} as the generalization error. The generalization error is introduced by substituting the true risk with the empirical risk, and it is the only component of the total error which depends on the sampling method. It is vital to derive the convergence rate of generalization error with respect to the sample size $m$.  
    
    Compared to the classical MC method, the QMC method provides more accurate approximations of the true risk as the sample size increases. Consequently, it is important to study whether QMC is a better choice for our problem. To this end, we focus on studying the convergence rates of generalization error \eqref{eq:generalization_error} under MC and QMC frameworks. We also numerically compare the efficiency of MC-based and QMC-based DBDP methods. These results are presented in Sections~\ref{sec:main} and~\ref{sec:numerical}.
	
	\section{The analysis of generalization error}\label{sec:main}
	
	\subsection{Problem specification}\label{sec:problem}
	We consider two basic types of nonlinear parabolic partial PDEs: the nonlinear heat PDEs and nonlinear Black-Scholes PDEs. For heat equations, the SDE is  
	\begin{equation}\label{eq:heat}
		\X_t = \eta + \overline{\mu} t+  \overline{\sigma} W_t, \;\;\; 0 \le t \le T ,
	\end{equation}
	where $\overline{\mu} \in \mathbb{R}^{d}, \overline{\sigma} \in \mathbb{R}^{d\times d}$. For Black-Scholes equations, the SDE is 
	\begin{equation}\label{eq:BS}
		\X_t = \eta + \int_0^t \mathrm{diag}(\X_s)\overline{\mu} \dd s+  \int_0^t  \mathrm{diag}(\X_s)\overline{\sigma} \dd W_s , \;\;\; 0 \le t \le T ,
	\end{equation}
	where $\mathrm{diag}(x)$ denotes the diagonal matrix with its diagonal elements equal to $x$. The solution of \eqref{eq:BS} is given by 
	\begin{align}\label{eq:BS_solu}
		\X_{t,i} = \eta^{(i)}\exp\left( \left( \mu_{i}-\frac{1}{2}\norm{\sigma_i}_{2}^{2}\right) t+\seq{\sigma_i,W_t}\right), \ \ 1\leq i \leq d,  
	\end{align}
	where $\X_{t,i}$, $\mu_i$ and $\eta^{(i)}$ denote the $i$-th components of $\X_t$, $\overline{\mu}$ and $\eta$,  respectively, and $\sigma_i$ denotes the $i$-th row of $\overline{\sigma}$. For simplicity, we take $\overline{\mu}=\mu\bm{1}_d$ and $\overline{\sigma}=\sigma\bm{I}_d$, where $\mu, \sigma \in \R$, $\bm{1}_d$ is the d-dimensional vector where all elements equal 1 and $\bm{I}_d$ is the d-dimensional identity matrix. All results in Section \ref{sec:main} can be extended to more general $\overline{\mu}$ and  $\overline{\sigma}$.
	
	For both nonlinear heat PDEs and nonlinear Black-Scholes PDEs, we can simulate $\X_{t_i}$ without using Euler-Maruyama method. More precisely, for the nonlinear heat PDEs, we use
	\begin{align}\label{eq:heat_formula}
		X_{t_{i+1}} &= X_{t_{i}} +\mu \Delta t_i\bm{1}_d+\sigma\Delta W_{t_i} \\
		X_{t_{i}} &= \eta + \mu t_i\bm{1}_d + \sigma W_{t_i}, \nonumber
	\end{align}
	and for nonlinear Black-Scholes PDEs, we use
	\begin{align}\label{eq:BS_formula}
		X_{t_{i+1},k} &= X_{t_{i},k}\exp\{(\mu-\sigma^2/2)\Delta t_i+\sigma\Delta W_{t_i,k} \} \\
		X_{t_{i},k} &= \eta_k \exp\{(\mu-\sigma^2/2)t_i+\sigma W_{t_i,k} \}, \nonumber
	\end{align}
	where $X_{t_{i},k}$, $\eta_k$ and $\Delta W_{t_i,k}$ denote the $k$-th components of $X_{t_{i}}$, $\eta$ and $\Delta W_{t_i}$, respectively. Note that $X_{t_i} = \X_{t_i}$. Therefore, Assumption~\ref{assump:2} holds.
	Since $\eta$, $W_{t_i}$ and $\Delta W_{t_i} $ are independent random variables, we can generate $(\eta$, $W_{t_i},\Delta W_{t_i})$ using a $3d$-dimensional standard normal random variable $W$, i.e.,
	\begin{align}\label{eq:Trans}
		\eta = (b-a)\Phi(W_{1:d}) + a \ \ , W_{t_i} = \sqrt{t_i}W_{(d+1):2d}, \ \ \Delta W_{t_i} = \sqrt{\Delta t_i}W_{(2d+1):3d},
	\end{align}
	where $\Phi(\cdot)$ is the cumulative distribution function of the the standard Gaussian distribution, acting on each component of $W_{1:d}$, and $x_{i:j}$ denotes $(x_i,\cdots,x_j)$. The loss function \eqref{agm:DBDP1} can be written as 
	\begin{equation}\label{eq:loss_function_simple}
		\mathcal{L}_i(\theta) = \E\left[ G_i(W;\theta)\right],
	\end{equation}
	where $G_i(\cdot;\theta)$ can be derived from~\eqref{agm:DBDP1} and~\eqref{eq:Trans}. The following lemma gives the upper bound and the lipschitz property of $G_i(W;\theta)$.
	\begin{lemma}\label{lem:G_property}
		Denote $\mS_1=(d,l_2,\cdots,l_{L},1)$ and $\mS_2=(d,l_2,\cdots,l_{L},d)$. Suppose that the neural networks $\mathcal{U}_i\in \mathcal{N}_{\mS_1,R}$ and $\mathcal{Z}_i\in \mathcal{N}_{\mS_2,R}$, $0\le i\le N-1$. For nonlinear heat PDEs and nonlinear Black-Scholes PDEs, there exist constants $C_1,C_2>0$ that depend only on $T,\mu,\sigma$, such that for every $i,\theta,\theta_1,\theta_{2}$, we have
		\begin{equation*}
			\begin{aligned}
				\abs{G_i(W;\theta)} &\leq C_1C_{fg} d^2B^2_{\mS_1,R}\exp\{ C_2\norm{W}_\infty\}, \\
				\abs{G_i(W;\theta_1)-G_i(W;\theta_2)}&\leq C_1C_{fg} d^2B_{\mS_1,R}C_{\mS_1,R} \norm{\theta_1-\theta_2}_\infty \exp\{ C_2\norm{W}_\infty\},
			\end{aligned}
		\end{equation*}
		where $C_{fg}=\max\{C_f^2,C_g^2,f^2(0),g^2(0)\}$, $C_{\mS_1,R}=|\mS_1|R^{D(\mS_1)-1}\norm{\mS_1}_\infty^{D(\mS_1)-1}$ and $B_{\mS_1,R}=2R\norm{\mS_1}_\infty$ (see Definition~\ref{defn:fnn} and Lemma~\ref{lem:property_FNN}).
	\end{lemma}
	
	\begin{proof}
		Denote $$H(X_{t_{i}},X_{t_{i+1}},\Delta W_{t_i};\theta) = \widehat \cU_{i+1} (X_{t_{i+1}}) - F(t_i,X_{t_i},\cU_i(X_{t_i};\theta),\cZ_i(X_{t_i};\theta),\Delta t_i,\Delta W_{t_i}).$$
		From \eqref{eq:F}, the lipschitz condition of functions $f,g$ and Lemma \ref{lem:property_FNN}, we have
		\begin{small}
			\begin{equation}
				\begin{aligned}
					&\abs{F(t, x, y, z, h, \Delta)} \leq\abs{y}+\abs{h}\abs{f(0)}+\abs{h}C_f\left(\sqrt{t} +\abs{y}+\norm{x}_2+\norm{z}_2\right) +\norm{z}_\infty\norm{\Delta}_1,\\
					&\abs{F(t, x, y_1, z_1, h, \Delta)-F(t, x, y_2, z_2, h, \Delta)} 
					\leq \abs{y_1-y_2}+\abs{h}C_f \abs{y_1-y_2}\\
					&\qquad \qquad \qquad \qquad \quad +\abs{h}C_f\sqrt{d}\norm{z_1-z_2}_\infty+\norm{z_1-z_2}_\infty\norm{\Delta}_1.
				\end{aligned}\nonumber
			\end{equation}
		\end{small}
		By applying Lemma \ref{lem:property_FNN}, we can prove that
		\begin{small}
			\begin{equation}\label{eq:H_property}
				\begin{aligned}
					&\big{|} H(X_{t_{i}},X_{t_{i+1}},\Delta W_{t_i};\theta)\big{|} 
					\leq \abs{\widehat \cU_{i+1} (X_{t_{i+1}})} + C_f\abs{\pi}\norm{X_{t_i}}_{2}+B_{\mS_1,R}\norm{\Delta W_{t_i}}_1\\
					&\qquad \qquad \qquad \qquad \qquad +(1+C_f\abs{\pi}+C_f\abs{\pi}\sqrt{d})B_{\mS_1,R}+\abs{f(0)}\abs{\pi}+\sqrt{T}C_f\abs{\pi}, \\
					&\big{|} H(X_{t_{i}},X_{t_{i+1}},\Delta W_{t_i};\theta_1)-
					H(X_{t_{i}},X_{t_{i+1}},\Delta W_{t_i};\theta_2)
					\big{|} \\
					&\leq  (1+C_f\abs{\pi}+C_f\abs{\pi}\sqrt{d}+\norm{\Delta W_{t_i}}_1)C_{\mS_1,R}\max\{\norm{X_{t_i}}_\infty,1\}\norm{\theta_1-\theta_2}_\infty.
				\end{aligned}
			\end{equation}
		\end{small}
		For each $i$, we have
		\begin{equation}
			\begin{aligned}\label{eq:hatU}
				\abs{\widehat \cU_{i+1} (X_{t_{i+1}})}
				&\leq \max\{B_{\mS_1,R},\abs{g(X_{t_{i+1}})}\} \\
				&\leq \max\{B_{\mS_1,R},\abs{g(0)} + C_g \norm{X_{t_{i+1}}}_2\}.
			\end{aligned}
		\end{equation}
		For the simulation scheme \eqref{eq:heat_formula} and \eqref{eq:Trans}, we can easily prove that
		\begin{equation}\label{eq:heat_simu}
			\begin{aligned}
				\norm{X_{t_{i}}}_{\infty} &\leq 1+\abs{\mu}T+\sigma\sqrt{T}\norm{ W}_{\infty}, \\
				\norm{X_{t_{i+1}}}_{\infty} &\leq 1+2\abs{\mu}T+2\sigma\sqrt{T}\norm{ W}_{\infty}, \\
				\norm{\Delta W_{t_i}}_{\infty} &\leq \abs{\pi}\norm{ W}_{\infty}.
			\end{aligned}
		\end{equation}
		For the simulation scheme \eqref{eq:BS_formula} and \eqref{eq:Trans}, we also have 
		\begin{equation}\label{eq:BS_simu}
			\begin{aligned}
				\norm{X_{t_{i}}}_{\infty} &\leq \exp\{\abs{\mu-\sigma^2/2}T+\sigma T \norm{ W}_{\infty} \},\\
				\norm{X_{t_{i+1}}}_{\infty} &\leq \exp\{2\abs{\mu-\sigma^2/2}T+2\sigma T \norm{ W}_{\infty} \}, \\
				\norm{\Delta W_{t_i}}_{\infty} &\leq \abs{\pi}\norm{ W}_{\infty}.
			\end{aligned}
		\end{equation}
		Combining \eqref{eq:H_property}, \eqref{eq:hatU}, \eqref{eq:heat_simu}, \eqref{eq:BS_simu} and the fact that
		$$G_i(W;\theta)=H^2(X_{t_{i}},X_{t_{i+1}},\Delta W_{t_i};\theta),$$
		the proof is completed after straightforward calculation.
	\end{proof}

\subsection{Generalization error for Monte Carlo methods}
	As shown in Section~\ref{sec:total_error}, the total error is bounded by the sum of three parts, namely the scheme error, the approximation error and the generalization error. In this section, we study the generalization error under the framework of MC. For simplicity, we drop the subscript $i$ in~\eqref{eq:loss_function_simple}.

	Assume that there exists the optimal parameter $\theta^{*}$ such that
	$$\theta^{*} = \arg\min_{\theta \in \Theta}\mathcal{L}(\theta).$$
	
	Since the expectation is usually intractable, we apply the ERM method to estimate $\theta^{*}$ by the minimizer $\theta^{(m)}$ of the empircal risk 
	\begin{equation}\label{eq:MC_empircalrisk}
		\mathcal{L}^{(m)}\left( \theta\right) = \frac{1}{m}\sum_{i = 1}^{m}G(\bm w_{i};\theta) 
	\end{equation}
	where $\left\{\bm w_{i}: 1\leq i \leq m \right\}$ are random samples satisfying $\bm w_{i} \sim W$. 
	
	We are interested in the generalization error defined by (see Theorem~\ref{thm:total_error} and Corollary~\ref{cor:error_decom} for details) 
	$$\mathcal{L}(\theta^{(m)})-\mathcal{L}(\theta^{*}).$$
	
	The next lemma establishes the upper bound of the generalization error.
	
	\begin{lemma}\label{lem:upperbound}
		For every $m\in \mathbb{N}$, we have
		\begin{eqnarray*}
			0\leq \mathcal{L}(\theta^{(m)})-\mathcal{L}(\theta^{*})&\leq& 2 \underset{\theta \in \Theta}{\mathrm{sup}} \abs{\mathcal{L}(\theta)-\mathcal{L}^{(m)}(\theta)}
		\end{eqnarray*}
	\end{lemma}
	\begin{proof}
		The lower bound follows by $\theta^{*} = \mathrm{argmin}_{\theta \in \Theta}\mathcal{L}(\theta)$. For the upper bounds, we have
		\begin{eqnarray*}
			\mathcal{L}(\theta^{(m)})-\mathcal{L}(\theta^{*})&=& \mathcal{L}(\theta^{(m)})-\mathcal{L}^{(m)}(\theta^{(m)}) + \mathcal{L}^{(m)}(\theta^{(m)})-\mathcal{L}^{(m)}(\theta^{*})\\
			&&+ \mathcal{L}^{(m)}(\theta^{*})-\mathcal{L}(\theta^{*})\\
			&\leq&\mathcal{L}(\theta^{(m)})-\mathcal{L}^{(m)}(\theta^{(m)}) +\mathcal{L}^{(m)}(\theta^{(m)})-\mathcal{L}^{(m)}(\theta^{*})\\
			&\leq&2 \underset{\theta \in \Theta}{\mathrm{sup}} \abs{\mathcal{L}(\theta)-\mathcal{L}^{(m)}(\theta)},
		\end{eqnarray*}
		where the second inequality follows by $\mathcal{L}^{(m)}(\theta^{(m)})\leq\mathcal{L}^{(m)}(\theta^{*})$.
	\end{proof}
	
	\subsubsection{The convergence rate of mean error for MC}
	Under the MC framework, $(\bm w_{i})_{i=1}^{m}$ in~\eqref{eq:MC_empircalrisk} are identical and independent standard Gaussian random variables. By Lemma~\ref{lem:upperbound}, the mean generalization error satisfies
	\begin{equation}\label{eq:mean_error_bound}
		\E\left[ \cL(\theta^{(m)}) - \cL(\theta^*)\right] \le 2 \E\left[ \sup_{\theta \in \Theta} \left| \frac{1}{m}\sum_{i=1}^m G(\bm w_i;\theta) - \E\left[ G(W;\theta)\right]\right|\right],
	\end{equation}
	where $W$ is a standard Gaussian random variable. As a result, it suffices to obtain the upper bound of the right hand side of~\eqref{eq:mean_error_bound}. To this end, we use the Rademacher complexity technique. 
	\begin{Definition}[Rademacher complexity]\label{defn:rc}
		Suppose that $( \bm x_{i})_{i=1}^{m}$ are i.i.d. random variables on $(\Omega,\mathcal{F},\mathbb{P})$ and take values in a metric space $S$. Suppose $\mathcal{G}$ is a class of real-valued borel-measurable function over $(S,\mathcal{B}(S))$ ($\mathcal{B}(S)$ is the Borel $\sigma$-field). The Rademacher complexity of a function class $\mathcal{G}$ with respect to the random sample $(\bm x_{i})_{i=1}^{m}$ is defined by
		\begin{equation}
			R_{m}(\mathcal{\mathcal{G}};( \bm x_{i})_{i=1}^{m}) = \E\left[ {\underset{g \in \mathcal{G}}{\mathrm{sup}}\frac{1}{n}\sum_{i=1}^m\xi_i g(\bm  x_{i})}\right] , \nonumber	
		\end{equation}
		where $\set{\xi_i}_{i=1}^{m}$ are i.i.d. Rademacher variables with discrete distribution $\mathbb{P}\{\xi_i=1\}=\mathbb{P}\{\xi_i=-1\}=1/2$. 
	\end{Definition}
 
	The next lemma provides an upper bound for the right hand side of \eqref{eq:mean_error_bound} in terms of the Rademacher complexity. We refer to \cite[Lemma 26.2]{shalev2014understanding} and \cite[Theorem 3.3] {mohri2018foundations} for a detailed proof.
	\begin{lemma}\label{lem:rcbasic}
		Under the setting of Definition \ref{defn:rc}, we have
		\begin{equation}
			\E\left[ {\underset{g \in \mathcal{G}}{\mathrm{sup}}\abs{\frac{1}{m}\sum_{i=1}^m g( \bm x_{i})-\E[g( \bm x)]}}\right]  \leq 2R_{m}(\mathcal{\mathcal{G}};( \bm x_{i})_{i=1}^{m}) .\nonumber
		\end{equation}
	\end{lemma}   
 
	The Rademacher complexity has the following upper bound. 
	\begin{lemma}\label{lem:Rand}
		Under the setting of Definition \ref{defn:rc}, suppose $\mathcal{G}$ is a parameterized class of functions $\left\lbrace \mathcal{F}(\cdot;\theta):\mathcal{S} \to \mathbb{R}\mid\theta \in \Theta\right\rbrace $, where  $\Theta \subseteq \mathbb{R}^{p}$ is the parametric space satisfying $\norm{\theta}_{\infty}\leq R$. Suppose that there exist constants $B_1$ and $B_2$ such that for all $\theta_{1},\theta_{2}\in \Theta$ it holds that
		\begin{eqnarray*}
			\underset{x \in S}{\mathrm{sup}}\abs{\mathcal{F}(x;\theta_{1})-\mathcal{F}(x;\theta_{2})} 
			&\leq& B_1\norm{\theta_{1}-\theta_{2}}_{\infty}, \\
			\underset{x \in S}{\mathrm{sup}}\abs{\mathcal{F}(x;\theta_{1})}
			&\leq& B_2.
		\end{eqnarray*}
		Then we have
		\begin{eqnarray*}
			R_{m}(\mathcal{G};( x_{i})_{i=1}^{m})&\leq& \frac{4}{\sqrt{m}}+\frac{6\sqrt{p}B_2}{\sqrt{m}}\sqrt{\log\left(2B_1R\sqrt{m} \right) }.
		\end{eqnarray*}
	\end{lemma}
	
	For a detailed proof of Lemma~\ref{lem:Rand}, we refer to \cite[Theorem 5.1]{jiao2024error}. By Lemma~\ref{lem:rcbasic}, the right hand side of~\eqref{eq:mean_error_bound} is bounded by the Rademacher complexity of function class $\{G(\cdot;\theta)\}_{\theta\in\Theta}$. However, we can not directly use Lemma~\ref{lem:Rand} to obtain the upper bound of this Rademacher complexity, since $G(\cdot;\theta)$ is not bounded. To address this issue, we use the truncation method to obtain the mean generalization error rate of MC.
	\begin{theorem}\label{thm:mc_SampleComplexity}
		In the cases of nonlinear heat PDEs and nonlinear Black-Sholes PDEs, suppose that we use the MC method to obtain the samples $\left\{\bm w_{i}: 1\leq i \leq m \right\}$ in \eqref{eq:MC_empircalrisk}, then we have
		\begin{align*}
			\E \left[  \mathcal{L}(\theta^{(m)})-\mathcal{L}(\theta^{*})\right] = Cm^{-1/2}\log(m),
		\end{align*}
		where $C$ is a polynomial of $(R,\|\mS_1\|_{\infty},D(\mS_1),C_{fg})$, independent of $m$.
	\end{theorem}
	\begin{proof}
		For every $K>0$, denote $G^{(K)}(x;\theta) := \bm 1_{\{\|x\|_{\infty}\le K\}}G(x;\theta)$ be the truncated version of $G(x;\theta)$. The truncated function class $\mathcal{G}^{(K)} :=\{G^{(K)}(x;\theta)\}_{\theta\in\Theta}$ satisfies the conditions of Lemma~\ref{lem:Rand}. It follows from~\eqref{eq:mean_error_bound} that
		\begin{align}
			\E \left[  \mathcal{L}(\theta^{(m)})-\mathcal{L}(\theta^{*})\right] 
			\leq& 2\underset{\theta \in \Theta}{\mathrm{sup}} \abs{\E\left[  \bm1_{\{\|W\|_{\infty} > K\}}G(W;\theta)\right]}\notag\\
			&+ 2\E \left[\underset{\theta \in \Theta}{\mathrm{sup}} \abs{\bm1_{\{\|W\|_{\infty} > K\}}G(W;\theta)}\right] +4R_{m}(\mathcal{G}^{(K)};( \bm w_{i})_{i=1}^{m})\notag\\
			\le& 4\E\left[\underset{\theta \in \Theta}{\mathrm{sup}} \abs{\bm1_{\{\|W\|_{\infty} > K\}}G(W;\theta)}\right] + 4R_{m}(\mathcal{G}^{(K)};( \bm w_{i})_{i=1}^{m}). \label{mideq_mc_mean_error:two_terms}
		\end{align}
		Observed by Lemma \ref{lem:G_property}, we have 
		\begin{align}
			\abs{G(W;\theta)} &\leq C_1C_{fg} d^2B^2_{\mS_1,R}\exp\{ C_2\|W\|_{\infty}\}\label{mideq:truncated_bound},\\
			\abs{G(W;\theta_1)-G(W;\theta_2)}&\leq C_1C_{fg} d^2B_{\mS_1,R}L_{\mS_1,R}\exp\{ C_2\|W\|_{\infty}\}\norm{\theta_1-\theta_2}_\infty.\label{mideq:lip}
		\end{align}
		For the first term on the right hand side of~\eqref{mideq_mc_mean_error:two_terms}, it follows from~\eqref{mideq:truncated_bound} and H\"older's inequality that 
		\begin{equation}\label{mideq:tail_bound}
			\begin{aligned}
				\E\left[\underset{\theta \in \Theta}{\mathrm{sup}} \abs{1_{\{\|W\|_{\infty} > K\}}G(W;\theta)}\right] &\le \sqrt{\E\left[\bm 1_{\{\|W\|_{\infty}>K\}} \right]}\sqrt{\E\left[ \sup_{\theta\in\Theta}|G(W;\theta)|^2\right]}\\
				&\le C_1C_{fg} d^3B^2_{\mS_1,R}\sqrt{d\E\left[\exp\{ C_2|W_1|\}\right]}\sqrt{\PP{\|W\|_{\infty} > K}}\\
				&\le \sqrt{2}C_1C_{fg} d^{4}B^2_{\mS_1,R}\sqrt{\E\left[\exp\{ C_2|W_1|\}\right]}\exp\{-K^2/4\},
			\end{aligned}
		\end{equation}
		where in the second inequality, we use the facts that $$ \E\left[\exp\{ C_2\|W\|_{\infty}\}\right]\le d\E\left[\exp\{ C_2|W_1|\}\right],$$ and 
		\begin{equation}\notag
			\PP{\norm{W}_\infty>K} \leq \sum_{i=1}^{d} \PP{\abs{W_i}>K} 
			\leq 2d\exp\{-K^2/2\}.
		\end{equation}
		For the second term on the right hand side of~\eqref{mideq_mc_mean_error:two_terms}, it follows from~\eqref{mideq:lip} and Lemma \ref{lem:Rand} that
		\begin{align}\label{mideq:radermacher_bound}
			R_{m}(\mathcal{G}^{(K)};( \bm w_{i})_{i=1}^{m}) 
			\leq \frac{4}{\sqrt{m}}+\frac{6\sqrt{|\mS_1|+|\mS_2|}B_2(K)}{\sqrt{m}}\sqrt{\log\left(2B_1(K)R\sqrt{m} \right) },
		\end{align}
		where $B_1(K) =C_1C_{fg} d^2B_{\mS_1,R}C_{\mS_1,R}\exp\{ C_2K\}$, 
		$B_2(K)= C_1C_{fg} d^2B^2_{\mS_1,R}\exp\{ C_2K\}.$
		
		We take $K=C_2^{-1}\log(m)$. The upper bound of~\eqref{mideq:tail_bound} becomes 
		\begin{equation}\label{mideq:tail_bound_bound}
			h_1\eexp{-(C_2^{-1}\log(m))^2/4} \le \eexp{C_2^2}h_1m^{-1/2},
		\end{equation} 
		where $h_1$ is a polynomial of $(d,R,\|\mS_1\|_{\infty},C_{fg})$ and we use the fact that $(C_2^{-1}\log(m))^2 + C_2^{2} \ge 2\log(m)$. The upper bound of~\eqref{mideq:radermacher_bound} becomes
		\begin{equation}\label{mideq:radermacher_bound_bound}
			h_2 m^{-1/2}\sqrt{\log(m\sqrt{m})} \le h_2 m^{-1/2}\log(m),
		\end{equation}
		where $h_2$ is a polynomial of $(d,R,\|\mS_1\|_{\infty},D(\mS_1),C_{fg})$. Combining~\eqref{mideq:tail_bound},~\eqref{mideq:tail_bound_bound},~\eqref{mideq:radermacher_bound} and~\eqref{mideq:radermacher_bound_bound}, we prove the desired result.
	\end{proof}
	
	If we use the MC method to obtain the samples for computing the empirical risk, Theorem \ref{thm:mc_SampleComplexity} demonstrate that the generalization error converges at the rate of $O(\log(m)m^{-1/2})$ with respect to the sample size $m$. This rate is consistent with those found in classical learning theory \cite{shalev2014understanding,mohri2018foundations}, where the output is usually assumed to be bounded. 
 
	\subsubsection{Break the curse of dimensionality}
        We denote $\mathcal{L}_i^* = L_i(\cU_{i+1},\cU_i^*,\cZ_i^*)$. For simplicity, we drop the subscript $ i$ in this section. In Theorem \ref{thm:mc_SampleComplexity}, we find that the implied constant is actually a polynomial of the constant $C_{fg}$ (see Lemma \ref{lem:G_property}) and the size of the neural networks, including the parameter bound $R$, the width $\|\mS_1\|_{\infty}$, the depth $D(\mS_1)$. This allows us to prove that the ERM method breaks the curse of dimensionality, provided that the approximation error in Corollary 
        \ref{cor:error_decom} satisfies certain assumptions.

        \begin{Assumption}\label{ass:approx_UZ}
            Let $\mS_1=(d,l_2,\cdots,l_{L},1)$, $\mS_2=(d,l_2,\cdots,l_{L},d)$. For each $i\in\{1,2,\cdots,N\}$, there exists a poynomial $p$ such that for for every $\delta>0$, there exists a neural network $\cU_\delta\in \mathcal{N}_{\mS_1,R}$ and $\cZ_\delta \in \mathcal{N}_{\mS_1,R}$ such that 
            \begin{equation}
                 \E\left[ \left| \cU_i^*(X_{t_i}) - \cU_\delta(X_{t_i})\right|^2\right] +\E\left[ \left| \cZ_i^*(X_{t_i}) - \cZ_\delta(X_{t_i})\right|^2\right] \leq \delta,
            \end{equation}
            and 
            \begin{equation}
                \max\{R,D(\mS_1),\norm{\mS_1}_\infty\}\leq p(d,\delta^{-1}).
            \end{equation}
        \end{Assumption}
        \begin{theorem}\label{thm:excessrisk}
		Under Assumption \ref{ass:approx_UZ}, suppose that the PDE \eqref{eq:PDE} is nonlinear heat PDEs or nonlinear Black-Sholes PDEs described in Section \ref{sec:problem}. Let $C_{fg}=\max\{C_f^2,C_g^2,f^2(0),g^2(0)\}$ ,$\mS_1=(d,l_2,\cdots,l_{L},1)$ and $\mS_2=(d,l_2,\cdots,l_{L},d)$. Then there exists a polynomial $h$ such that for $m$ satisfying
		\begin{equation*}
			m\geq h\left( \delta^{-1},d,C_{fg}\right)   ,
		\end{equation*}
		we have
		\begin{equation}\label{eq:expecrbound}
			\E[\mathcal{L}(\theta^{(m)})-\mathcal{L}^{*}] \leq \delta .
		\end{equation}
		Moreover, the coefficients of $h$ only depends on the constants $C_1,C_2$ (see Lemma~\ref{lem:G_property}).
	\end{theorem}

        \begin{proof}
            This is a direct consequence of Corollary \ref{cor:error_decom} and Theorem \ref{thm:mc_SampleComplexity}.
         \end{proof}
         
	Studying the probability inequality of the excess risk $ \mathcal{L}(\theta^{(m)})-\mathcal{L}^*$ is an important topic in  the analysis of the ERM method, see \cite{shalev2014understanding,mohri2018foundations,berner2020analysis,jiao2024error,zhang2004statistical}. In some contexts, the problem is for $\delta, \rho > 0$, determining the $m$ such that 
	\begin{equation}\label{eq:probablity_ineq}
		\PP{\mathcal{L}(\theta^{(m)})-\mathcal{L}^{*} \leq \delta} \geq 1-\rho.
	\end{equation}
    It follows from Markov inequality that for $\delta > 0$,\
    \begin{equation}\notag
        \PP{\cL(\theta^{(m)}) - \mathcal{L}^{*} \le \delta } \ge 1- \frac{\E\left[ \cL(\theta^{(m)}) - \mathcal{L}^{*}\right]}{\delta}.
    \end{equation}
    Thus we have the natural extension of Theorem \ref{thm:excessrisk}.
    \begin{Corollary}\label{cor:total_error_2}
        Under the setting of Theorem \ref{thm:excessrisk}. There exists a polynomial $q$ such that for $m$ satisfying
        \begin{equation*}
            m\geq q\left( \delta^{-1},\rho^{-1},d,C_{fg}\right)   ,
        \end{equation*}
        we have
        \begin{equation}
            \PP{\mathcal{L}(\theta^{(m)})-\mathcal{L}^{*} \leq \delta} \geq 1-\rho .
        \end{equation}
        Moreover, the coefficients of $q$ only depends on the constants $C_1,C_2$ (see Lemma~\ref{lem:G_property}).
    \end{Corollary}
	
    In fact, Theorem \ref{thm:excessrisk} and Corollary \ref{cor:total_error_2} show that, if the functions $\cU^{*}$ and $\cZ^{*}$ can be approximated by the neural networks without the curse of dimensionality, then the required sample size $m$ to satisfy the probabilistic inequality grows polynomially, which means the total error of the ERM principle over the neural networks overcomes the curse of dimensionality as dimension $d$ increases. 

    There are numerous theoretical results on the approximation properties of the neural networks. For universal approximation results, we refer the readers to \cite{elbrachter2021deep,guhring2020error,siegel2020approximation,hornik1993some,hornik1991approximation}. For the approximation of solutions to PDEs, various studies have demonstrated that neural networks can overcome the curse of dimensionality, see \cite{grohs2023proof,jentzen2021proof,hutzenthaler2020proof,hutzenthaler2020overcoming,darbon2016algorithms,hutzenthaler2023overcoming}. Verifying the Assumption \ref{ass:approx_UZ} is challenging in general cases since the target functions are represented by expectations. Additionally, it remains an open problem  to prove that the neural networks can overcome the curse of dimensionality when approximating the solutions to the general forms of nonlinear heat equations and Black-Scholes equations. 
	
    \subsection{Generalization error for quasi-Monte Carlo}
        In this section, we consider the QMC method and study the corresponding generalization error rate.  

    \subsubsection{Quasi-Monte Carlo methods}\label{sec:qmc}
    QMC~\cite{Niederreiter1992} is an efficient quadrature method to numerically solve the integral problems on the unit cube $ [0,1]^d$. While the MC method, relying on random sampling, achieves a convergence rate of $O(m^{-1/2})$ with $m$ quadrature points, QMC leverages low discrepancy sequences, such as the Halton sequence, Faure sequence and Sobol' sequence~\cite{caflisch1997,Niederreiter1992,owen2013}, which facilitates a faster convergence rate of $O(m^{-1+\varepsilon})$ for an arbitrarily small $\varepsilon > 0$. This enhanced convergence rate of QMC is grounded in the Koksma-Hlawka inequality~\cite{Hlawka1972}, linking the error to both the variation of the integrand in the sense of Hardy and Krause and the uniformity of sample points as measured by star discrepancy. More precisely, the Koksma-Hlawka inequality states that
    \begin{equation}
        \left| \frac{1}{m}\sum_{i=1}^m f(\y_i) - \int_{[0,1]^d} f(\y)\dd \y\right| \le V_{\mathrm{HK}}(f)D_m^*(\{\y_1,\dots,\y_m\}),
    \end{equation}
    where $V_{\mathrm{HK}}(f)$ is the variation of $f$ in the sense of Hardy and Krause and $D_m^*$ is the star discrepancy of point set $\{\y_1,\dots,\y_m\}$ with the magnitude $O(m^{-1}(\log m)^{d-1})$ for low discrepancy point sets (see~\cite{Niederreiter1992} for more details). The convergence rate of QMC is validated if the function $f$ is of bounded variation in the sense of Hardy and Krause (BVHK).

    The randomized quasi-Monte Carlo (RQMC) method, which includes random shifts~\cite{kuo2006a} and scrambling~\cite{owen1997b}, merges the benefits of randomization with low discrepancy sequences to obtain unbiased estimations. If the integrand is bounded and smooth, the nested scrambled method can achieve a faster convergence rate of $O(n^{-3/2+\varepsilon})$, as detailed by~\cite{owen2008}. Recently, Ouyang et al.~\cite{ouyang2023} have demonstrated that this high convergence rate can also be maintained for unbounded and smooth integrands by employing appropriate importance sampling methods.
    
    \subsubsection{The convergence rate of mean error for QMC}
    Combining QMC methods, the empirical risk is defined by
	\begin{equation}\label{eq:L_qmc}
		\widehat{\mathcal{L}}^{(m)}(\theta) := \frac{1}{m}\sum_{i=1}^m G(\Phi^{-1}(\y_i);\theta),
	\end{equation}
	where $\{\y_1,\dots,\y_m\}$ is a low discrepancy point set and $\Phi^{-1}$ is the inverse cumulative distribution function of standard Gaussian distribution, acting on each component of $\y_i$. In this paper, we consider the RQMC method with the scrambling approach (see Owen~\cite{owen1997b} for more details).  Define
	\begin{equation}
		\widehat \theta^{(m)} := \arg\min_{\theta\in \Theta} \widehat{\mathcal{L}}^{(m)}(\theta).
	\end{equation}

By Lemma~\ref{lem:upperbound}, the mean generalization error satisfies
\begin{equation}\label{eq:qmc_gen_error_bound}
     \begin{aligned}
	\E\left[ \cL(\widehat\theta^{(m)}) - \cL(\theta^*)\right] &\le 2\E\left[ \sup_{\theta \in \Theta} \Big| \cL(\theta) - \hcL_m(\theta)\Big|\right] \\
    &= 2\E\left[ \sup_{\theta \in \Theta}\Big| \E\left[ G(W;\theta)\right] - \sum_{i=1}^mG(\Phi^{-1}(\y_i);\theta)\Big|\right].
\end{aligned}
\end{equation}
Since $ G(\Phi^{-1}(\cdot);\theta)$ is an unbounded function (it can not be of BVHK), the Koksma-Hlawka inequality can not be directly used to provide a convergence rate for the right hand side of~\eqref{eq:qmc_gen_error_bound}. To address this issue, Xiao et al.~\cite{xiao2023analysis} applied the results in Owen~\cite{owen2006a} and obtained the convergence rate for the case of linear Kolmogorov PDEs.  Utilizing Owen~\cite{owen2006a} method would become quite cumbersome for our case, given that the function $ G$ is more complex than those in Xiao et al~\cite{xiao2023analysis}. Here, we recommend employing the results in Ouyang et al.~\cite{ouyang2023}, where we only need to examine the growth conditions of $ G$. We introduce the following function class from~\cite{ouyang2023}.

\begin{Definition}
		For $ h:\R^{3d} \rightarrow \R$, define the exponential growth function class as 
		\begin{equation}
			\mathcal{G}_e(A,B) := \left\{ h\in \mathcal{S}^{3d}(\R^{3d}) : \sup_{\uu \subset 1:3d} |\partial^{\uu}h(\x)|\le Be^{A|\x|}\right\},\notag
		\end{equation}
		where $ \uu \subset 1:3d = \left\{ 1,\dots,3d\right\}$ and $ \partial^{\uu} := \prod_{i\in \uu} \partial/\partial x_i$, $ \mathcal{S}^{3d}(\R^{3d})$ is the function class where the functions have continuous mixed partial derivatives of order $ 1$.
	\end{Definition}
 The following lemma can be obtained from the proof of Corollary 4.10 in Ouyang et al.~\cite{ouyang2023}, which gives the convergence rate of QMC methods regarding the function class $ \mathcal{G}_e(A,B)$.
	\begin{lemma}\label{lem:qmc_converge}
		Let $ \left\{ \y_1,\dots,\y_m\right\}$ be a RQMC point set such that each $ \y_j $ is uniformly distributed on $[0,1]^{3d}$ and 
		\begin{equation}
			D^*_m(\left\{ \y_1,\dots,\y_m\right\}) \le C\frac{(\log m)^{d-1}}{m} \quad a.s.,
		\end{equation}
		where $ D^*_m$ is the star discrepancy and $ C$ is a constant independent of $ m$. Then we have 
		\begin{equation}
			\E\left[ \sup_{h \in \mathcal{G}_e(A,B)} \Big|\sum_{i =1}^m h(\Phi^{-1}(\y_i)) - \E\left[ h(W)\right]\Big|\right] = O(m^{-1+\varepsilon}),
		\end{equation}
		where $ W$ is a standard $ 3d$-dimensional Gaussian random variable.
	\end{lemma}
 
 Noting that the nested scrambled Sobol' point set satisfies the conditions in Lemma~\ref{lem:qmc_converge}, if we use this QMC point set $ \{\y_i\}$ in~\eqref{eq:L_qmc}, then by Lemma~\ref{lem:qmc_converge} and~\eqref{eq:qmc_gen_error_bound}, it suffices to verify that the function $ G(\cdot;\theta)$ is in the class $ \mathcal{G}_e(A,B)$ for some constant $ A,B$. The main result of this section is the following theorem.
	\begin{theorem}\label{thm:qmc}
		If we use nested scrambled Sobol' point set in~\eqref{eq:L_qmc}, then for the nonlinear heat PDEs and the nonlinear Black-Scholes PDEs, we have 
		\begin{equation}
			\E\left[ \cL(\widehat \theta^{(m)}) - \cL(\theta^*)\right] = O(m^{-1+\varepsilon})
		\end{equation}
		for every $\varepsilon > 0$.
	\end{theorem}

	\begin{proof}
		By the same method used in Lemma~\ref{lem:upperbound}, we have
		\begin{equation}
			0\le \cL(\widehat\theta^{(m)}) - \cL(\theta^*) \le 2\sup_{\theta \in \Theta} \left| \cL(\theta) - \hcL_m(\theta)\right|.
		\end{equation}
		Thus,
		\begin{equation}\label{mideq:qmc11}
			\begin{aligned}
				\E\left[ \cL(\widehat\theta^{(m)}) - \cL(\theta^*)\right] &\le 2\E\left[ \sup_{\theta \in \Theta} \Big| \cL(\theta) - \hcL_m(\theta)\Big|\right] \\
				&= 2\E\left[ \sup_{\theta \in \Theta}\Big| \E\left[ G(W;\theta)\right] - \sum_{i=1}^mG(\Phi^{-1}(\y_i);\theta)\Big|\right].
			\end{aligned}
		\end{equation}
		We consider the growth of $ G(\cdot;\theta)$. 
		Note that $ G(\cdot;\theta)$ is composed through addition, multiplication and composition of some base functions (such as quadratic functions, linear functions, activation functions), where the derivatives of any order of these functions are dominated by a polynomial. By the results in Appendix of Ouyang et al.~\cite{ouyang2023}, there are constants $ A,B$, such that for any $ \theta \in \Theta $,
		\begin{equation}
			G(\cdot;\theta) \in \mathcal{G}_e(A,B). 
		\end{equation}
		Note that the nested scrambled Sobol' point set satisfies the conditions in Lemma~\ref{lem:qmc_converge}. By Lemma~\ref{lem:qmc_converge} and combining~\eqref{mideq:qmc11}, we have 
		\begin{equation}\notag
			\begin{aligned}
				\E\left[ \cL(\widehat\theta^{(m)}) - \cL(\theta^*)\right] &\le 2\E\left[ \sup_{\theta \in \Theta}\Big| \E\left[ G(W;\theta)\right] - \sum_{i=1}^mG(\Phi^{-1}(\y_i);\theta)\Big|\right] \\
				&\le 2\E\left[ \sup_{h \in \mathcal{G}_e(A,B)} \Big|\sum_{i =1}^m h(\Phi^{-1}(\y_i)) - \E\left[ h(W)\right]\Big|\right] = O(m^{-1+\varepsilon}).
			\end{aligned}
		\end{equation}
		This proves the desired result.
	\end{proof}

        Theorem~\ref{thm:mc_SampleComplexity} and~\ref{thm:qmc} establish the convergence rates of the generalization error under the MC and QMC methods. The results that the QMC method possesses a higher convergence rate $O(m^{-1+\varepsilon})$, indicating that, with sufficiently thorough training (achieving optimal parameters for each), the QMC method can achieve a smaller generalization error compared to the MC method. According to Theorem~\ref{thm:total_error} and Corollary~\ref{cor:error_decom} (which deals with the decomposition of the total error), the QMC method potentially enables the final trained solution to have a smaller error.  In the next section, we validate our theoretical results with several numerical examples.
	
	\section{Numerical experiments}\label{sec:numerical}
	In this section, we conduct numerical experiments to compare the performance of the DBDP scheme using different sample methods, specifically MC and RQMC methods. To ensure a fair comparison, we use the same training settings for different sampling methods, and to minimize the approximation error and highlight the impact of the generalization error, we choose a network larger than that in ~\cite{germain2022approximation}. For each step, we use a feedforward neural network with a depth of 4 and a width of $d+20$. The activation function is the Tanh function \eqref{eq:tanh}, and the network is initialized by means of Xavier initialization~\cite{glorot2010understanding}. The networks are trained using the Adam optimizer~\cite{kingma2014adam} with default weight decay $0.01$. For the first step near the terminal condition, the initial learning rate of is $0.01$ and is halved every $5000$ iterations, with a total of $50000$ iterations. For the subsequent steps, the initial learning rate of is $0.001$ and is halved every $500$ iterations, with a total of $5000$ iterations. We also add batch normalization layers~\cite{ioffe2015batch} to enhance the efficiency and robustness of the network training. 
	
	We measure the performance of the obtained solution by the relative $L_2$ error defined as
	\begin{equation}\label{eq:l2_error}
		\sqrt{\frac{\E\left[ \big|\widehat{\cU}(\eta)-u(0,\eta)\big|^2\right] }{\E[|u(0,\eta)|^2] }},
	\end{equation}
	where $\widehat{\cU}$ is the trained neural network, $u(0,\cdot)$ is the exact solution of the nonlinear parabolic PDEs \eqref{eq:PDE}, $\eta$ is uniformly distributed over target domain $[a,b]^d$. Since the relative $L_2$ error has no closed form, we approximated it by using the MC method as follows
	\begin{equation}\label{eq:relaErr}
		\sqrt{\frac{\sum_{i = 1}^{m} [\widehat{\cU}(\eta_i)-u(0,\eta_i)]^2 }{\sum_{i = 1}^{m} u^2(0,\eta_i)}},
	\end{equation}
	where $(\eta_i)_{i=1}^m$ are i.i.d. uniform samples on $[a,b]^d$. In the following experiments, we take the sample size $m=2^{16}$. 
 
        Note that unlike the solution obtained by training only at a single point in~\cite{germain2022approximation,han2018solving,hure2020}, we train solutions over an entire region. To demonstrate the efficiency of the training results at single points, in addition to presenting the $ L_2$ error of the results, we also provide histograms of the single-point errors.

	\subsection{Nonlinear heat equations}\label{sec:semi-heat}
	The first example is a nonlinear heat PDE  from~\cite{germain2022approximation}. We take the drift function $M(t,x) = \mu\bm{1}_d$, the diffusion function $\Sigma(t,x) = \sigma I_d$, the terminal condition $g(x) = \cos(\overline{x})$ with $\overline{x} = \sum_{i=1}^d x_i$, and the nonlinear term 
	\begin{equation}\label{eq_example:heat1}
		\begin{aligned}
			f(t,x,y,z) =& \left( \cos(\overline{x}) + \mu d \sin(\overline{x})\right)e^{\frac{T-t}{2}} - \frac{1}{2}\left( \sigma\sqrt{d}\sin(\overline{x})\cos(\overline{x})e^{T-t}\right)^2 \\
			&+\frac{1}{2d}(y(\bm1_d\cdot z))^2.
		\end{aligned}
	\end{equation}
	The solution of this PDE is $u(t,x) = \cos(\overline{x})e^{\frac{T-t}{2}}$. In our numerical experiment, we choose $\mu = 0.2/d$, $\sigma = 1/\sqrt{d}$ and the region $[a,b]^d=[-0.5,0.5]^d$. 
	
	Table \ref{tab:Heat_combined} presents the mean and standard deviation of the relative $L^{2}$ error over 8 independent runs for solving the nonlinear heat equation~\eqref{eq_example:heat1}, for dimensions 20 and 50, respectively. The table demonstrates that the RQMC-based algorithm outperforms the MC-based one. When the batch size is $2^{14}$, the mean error of RQMC-based solution is only $1/2$ of the MC-based solution. The standard deviations presented in Table~\ref{tab:Heat_combined} show that the variance of RQMC-based solutions is much smaller than the MC ones with the variance reduction factor (the ratio of MC variance to RQMC variance) varying from $4$ to $30$.
	
\begin{figure}[H]
		\centering
		\subfloat[Best trained neural networks]{
			\includegraphics[width=0.45\linewidth]{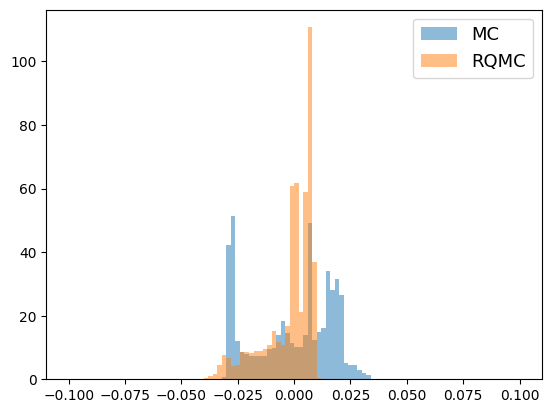}}
		\hspace{2em}
		\subfloat[Worst trained neural networks]{
			\includegraphics[width=0.45\linewidth]{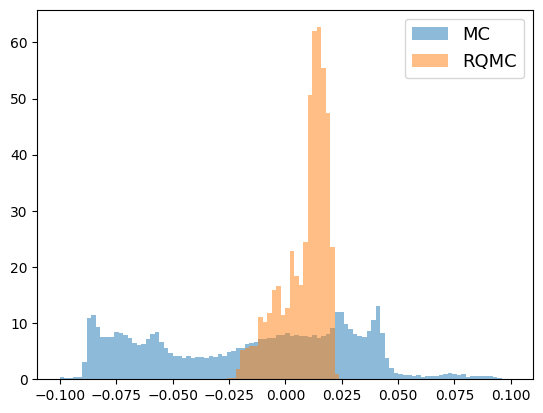}}
		\caption{The density histograms of the pointwise absolute errors for the best and worst trained neural networks in solving the nonlinear heat equation. The batch size is 16,384 and the  dimension is 50.}
		\label{fig:ErrorHist_heat}
	\end{figure}

	Figure \ref{fig:ErrorHist_heat} presents the density histogram of pointwise errors between the trained networks and the exact solution across the sample $(\eta_i)_{i=1}^{65536}$. The batch size is set to $2^{14}$, and the neural networks with the largest and smallest $L_2$ errors are selected from 8 independent runs. For the MC method, the smallest $L_2$ error is $1.65 \times 10^{-2}$ and the largest $L_2$ error is $4.12 \times 10^{-2}$. For the RQMC method, the smallest $L_2$ error is $9.32 \times 10^{-3}$ and the largest $L_2$ error is $1.14 \times 10^{-2}$. The figure demonstrates that the RQMC-based solution achieves smaller pointwise errors across the majority of samples. Figure \ref{fig:ErrorHist_heat} only displays the errors within the range of $[-0.1,0.1]$. For the best trained networks using MC methods, 143 out of 65,536 errors fall outside this range , while for the RQMC method, 0 out of 65,536 errors fall outside this range. This indicates that the solution based on RQMC approximates the exact solution with smaller error and smaller variance.

\begin{table}[H]
    \centering
    \begin{tabular}{|c|c|c|c|c|}
        \toprule
        Dimension & Batch Size & Method & Mean Error & Standard Deviation \\
        \midrule
        \multirow{6}{*}{20} & \multirow{2}{*}{$2^{12}$} & MC & $3.4 \times 10^{-2}$ & $6.8 \times 10^{-3}$ \\
        & & RQMC & $\bm{3.2 \times 10^{-2}}$ & $\bm{8.7 \times 10^{-4}}$ \\
        \cmidrule{2-5}
        & \multirow{2}{*}{$2^{14}$} & MC & $1.5 \times 10^{-2}$ & $6.0 \times 10^{-3}$ \\
        & & RQMC & $\bm{6.9 \times 10^{-3}}$ & $\bm{1.3 \times 10^{-3}}$ \\
        \cmidrule{2-5}
        & \multirow{2}{*}{$2^{16}$} & MC & $7.2 \times 10^{-3}$ & $1.0 \times 10^{-3}$ \\
        & & RQMC & $\bm{5.8 \times 10^{-3}}$ & $\bm{3.6 \times 10^{-4}}$ \\
        \midrule
        \multirow{6}{*}{50} & \multirow{2}{*}{$2^{12}$} & MC & $5.6 \times 10^{-2}$ & $1.8 \times 10^{-2}$ \\
        & & RQMC & $\bm{3.5 \times 10^{-2}}$ & $\bm{3.7 \times 10^{-3}}$ \\
        \cmidrule{2-5}
        & \multirow{2}{*}{$2^{14}$} & MC & $2.7 \times 10^{-2}$ & $7.3 \times 10^{-3}$ \\
        & & RQMC & $\bm{1.0 \times 10^{-2}}$ & $\bm{6.4 \times 10^{-4}}$ \\
        \cmidrule{2-5}
        & \multirow{2}{*}{$2^{16}$} & MC & $1.2 \times 10^{-2}$ & $2.8 \times 10^{-3}$ \\
        & & RQMC & $\bm{9.3 \times 10^{-3}}$ & $\bm{2.6 \times 10^{-3}}$ \\
        \bottomrule
    \end{tabular}
    \caption{The mean and standard deviation of the relative $L^{2}$ error over 8 independent runs for solving the nonlinear heat equation.}
    \label{tab:Heat_combined}
\end{table}

	\subsection{Hamilton--Jacobi--Bellman (HJB) equations}
	The second example is an HJB equation from \cite{beck2021Deep,han2018solving}. This specific HJB equation is given by
	\begin{equation}\label{eq_example:hjb}
		\begin{aligned}
			\frac{\partial u}{\partial t}(t,x) + (\Delta_xu)(t,x) & = \norm{(\nabla_xu)(t,x)}_2^2 , \;\;\;\;\;\; \text{ for  } (t,x)\in [ 0,T)\times\R^d, \\
			u(T,x) &=\norm{x}_2^{1/2}, \;\;\;\;\;  \mbox{ for } x \in \R^d.
		\end{aligned}
	\end{equation}
	This HJB equation can be transformed into a linear kolmogorov PDE through the Cole-Hopf transformation~\cite{cole1951quasi}. The solution can be expressed as 
    \begin{equation}
        u(t,x) = -\ln \left(\E\left[ \exp\left(-\|x+\sqrt{2}W_{T-t}\|_2^{1/2}\right)\right]\right).
    \end{equation}
    Therefore, to approximate the relative $L_2$ error~\eqref{eq:relaErr}, we use MC method to approximate the exact solution pointwise on samples $(\eta_i)_{i=1}^m$. In our numerical experiment, we choose the region $[a,b]^d=[0,1]^d$. 
	
\begin{figure}[H]
		\centering
		\subfloat[Best trained neural networks]{
			\includegraphics[width=0.45\linewidth]{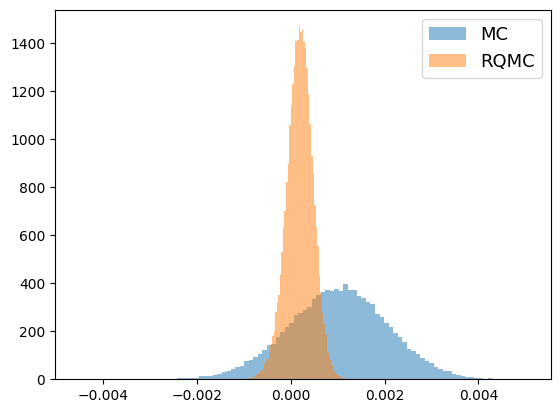}}
		\hspace{2em}
		\subfloat[Worst trained neural networks]{
			\includegraphics[width=0.45\linewidth]{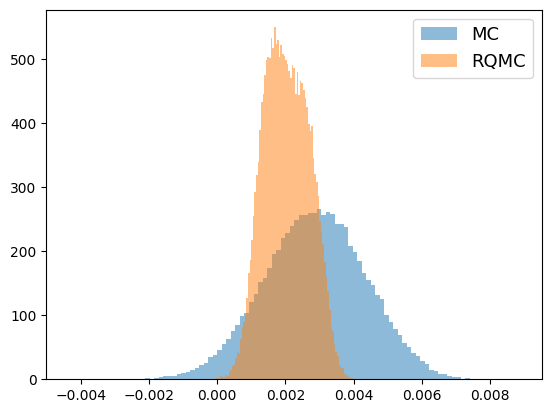}}
		\caption{The density histograms of the pointwise absolute errors for the best and worst trained neural networks in solving the HJB equation. The batch size is 16,384 and the dimension is 50.}
		\label{fig:ErrorHist_HJB}
	\end{figure}

	For the HJB equation~\eqref{eq_example:hjb}, Table \ref{tab:HJB_combined} presents the mean  and standard deviation of relative $L_2$ errors over 8 independent runs for dimensions 20 and  50, respectively. The table reveals that RQMC results are superior to MC, with a maximum error reduction to $1/4$ of MC error (e.g., when the dimension is 20 and the batch size is $2^{12}$). Figure \ref{fig:ErrorHist_HJB} presents the density histogram of pointwise errors between the trained networks and the exact solution. The batch size is set to $2^{14}$, and the neural networks with the largest and smallest $L_2$ errors are selected from 8 independent runs. For the MC method, the smallest $L_2$ error is $4.47 \times 10^{-4}$ and the largest $L_2$ error is $9.89 \times 10^{-4}$. For the RQMC method, the smallest $L_2$ error is $1.05 \times 10^{-4}$ and the largest $L_2$ error is $6.54 \times 10^{-4}$. Figure \ref{fig:ErrorHist_HJB} indicates that the distributions of errors roughly follow a normal distribution, with RQMC exhibiting smaller variance and a mean closer to 0. This suggests that the solutions trained using RQMC are closer to the true solutions with smaller errors, resulting in more stable outcomes. 

\begin{table}[H]
    \centering
    \begin{tabular}{|c|c|c|c|c|}
        \toprule
        Dimension & Batch Size & Methods & Mean Error & Standard Deviation \\
        \midrule
        \multirow{6}{*}{20} & \multirow{2}{*}{$2^{12}$} & MC & $1.7 \times 10^{-3}$ & $8.5 \times 10^{-4}$ \\
        & & RQMC & $\bm{4.4 \times 10^{-4}}$ & $\bm{2.1 \times 10^{-4}}$ \\
        \cmidrule{2-5}
        & \multirow{2}{*}{$2^{14}$} & MC & $1.2 \times 10^{-3}$ & $3.1 \times 10^{-4}$ \\
        & & RQMC & $\bm{3.9 \times 10^{-4}}$ & $\bm{2.0 \times 10^{-4}}$ \\
        \cmidrule{2-5}
        & \multirow{2}{*}{$2^{16}$} & MC & $5.4 \times 10^{-4}$ & $3.3 \times 10^{-4}$ \\
        & & RQMC & $\bm{1.9 \times 10^{-4}}$ & $\bm{9.1 \times 10^{-5}}$ \\
        \midrule
        \multirow{6}{*}{50} & \multirow{2}{*}{$2^{12}$} & MC & $1.3 \times 10^{-3}$ & $4.4 \times 10^{-4}$ \\
        & & RQMC & $\bm{2.8 \times 10^{-4}}$ & $\bm{5.3 \times 10^{-5}}$ \\
        \cmidrule{2-5}
        & \multirow{2}{*}{$2^{14}$} & MC & $6.7 \times 10^{-4}$ & $1.8 \times 10^{-4}$ \\
        & & RQMC & $\bm{2.9 \times 10^{-4}}$ & $\bm{1.5 \times 10^{-4}}$ \\
        \cmidrule{2-5}
        & \multirow{2}{*}{$2^{16}$} & MC & $3.6 \times 10^{-4}$ & $1.3 \times 10^{-4}$ \\
        & & RQMC & $\bm{1.7 \times 10^{-4}}$ & $\bm{8.3 \times 10^{-5}}$ \\
        \bottomrule
    \end{tabular}
    \caption{The mean and standard deviation of the relative $L^{2}$ error over 8 independent runs for solving the HJB equation.}
    \label{tab:HJB_combined}
\end{table}

	\subsection{Nonlinear Black-Scholes equations}
	
	Based on the nonlinear heat equations in Section \ref{sec:semi-heat}, we consider the corresponding nonlinear Black-Scholes PDE. We take the drift function $M(t,x) = \mu x$, the diffusion function $\Sigma(t,x) = \sigma \diag{x}$, and the terminal condition $g(x) = \cos(\overline{x})$ with $\overline{x} = \sum_{i=1}^d x_i$.
	
	For the first case, we take
	\begin{equation}\label{eq:f_BS_1}
		\begin{aligned}
			f(t,x,y,z) 
			=& \Big(\frac{1}{2}\cos(\overline{x}) + \mu\overline{x}\sin(\overline{x}) + \sigma^2\big(\sum_{i=1}^d x_i^2\big)\cos(\overline{x})\Big)e^{\frac{T-t}{2}}\\
			&- \frac{1}{2d}\left( \sigma\overline{x}\cos(\overline{x})\sin(\overline{x})e^{T-t}\right)^2 + \frac{1}{2d}(y(\bm1_d\cdot z))^2.
		\end{aligned}
	\end{equation}
	For \eqref{eq:f_BS_1}, the exact solution is $u(t,x) = \cos(\overline{x})e^{\frac{T-t}{2}}$.
	
	For the second case, we take a different $f$ from~\eqref{eq:f_BS_1}, which is
	\begin{equation}\label{eq:f_BS_2}
		\begin{aligned}
			f(t,x,y,z) =& \left( \frac{1}{2}\Phi(\overline{x}) - \mu \overline{x}\varphi(\overline{x}) + \sigma^2 \sum_{i=1}^d x_i^2 \varphi(\overline{x})\right)e^{\frac{T-t}{2}}\\
			& - \frac{1}{2d}\left( \sigma \overline{x}\varphi(\overline{x})\Phi(\overline{x})e^{T-t}\right)^2 + \frac{1}{2d}(y(\bm1_d\cdot z))^2,
		\end{aligned}
	\end{equation}
	where $\Phi$ is the cumulative distribution function of the one-dimensional standard Gaussian distribution and $\varphi$ is the corresponding probability density function. We take the terminal function $g = \Phi(\overline{x})$. The solution of this nonlinear Black-Scholes PDE is $u(t,x) = \Phi(\overline{x})e^{\frac{T-t}{2}}$.
	
	For the nonlinear term $\eqref{eq:f_BS_1}$, we take the region $[a,b]^d=[-0.5,0.5]^d$. Table \ref{tab:bs1_combined} presents the mean and standard deviation of the relative $L^{2}$ error over 8 independent runs, for dimensions 20 and  50, respectively. It demonstrates that the RQMC-based algorithm outperforms the MC-based one. Figure \ref{fig:ErrorHist_HJB} presents the density histogram of pointwise errors between the trained networks and the exact solution. The batch size is set to $2^{14}$, and the neural networks with the largest and smallest $L_2$ errors are selected from 8 independent runs. For the MC method, the smallest $L_2$ error is $1.74 \times 10^{-2}$ and the largest $L_2$ error is $4.58 \times 10^{-2}$. For the RQMC method, the smallest $L_2$ error is $9.01 \times 10^{-3}$ and the largest $L_2$ error is $1.64 \times 10^{-2}$. The figure demonstrates that the RQMC-based algorithm achieves smaller errors across the majority of samples. The histogram only displays the errors within the range of $[-0.1,0.1]$. For the best trained networks using MC methods, 138 out of 16,384 errors fall outside this range, while for the RQMC method, only 2 out of 16,384 errors fall outside this range.
	
\begin{figure}[H]
		\centering
		\subfloat[Best trained neural networks]{
			\includegraphics[width=0.45\linewidth]{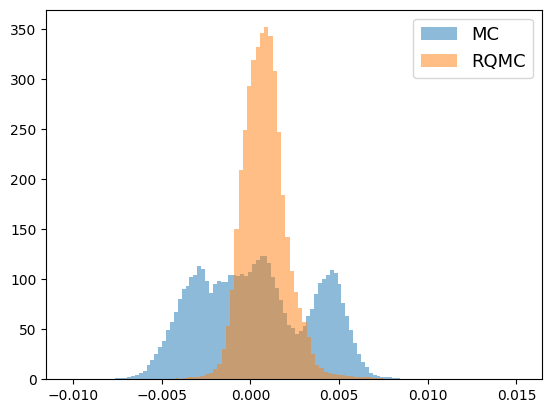}}
		\hspace{2em}
		\subfloat[Worst trained neural networks]{
			\includegraphics[width=0.45\linewidth]{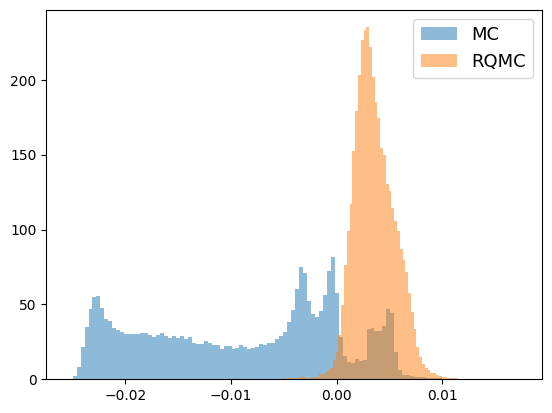}}
		\caption{The density histograms of the pointwise absolute errors for the best and worst trained neural networks in solving the  nonlinear Black-Scholes equation with the nonlinear term \eqref{eq:f_BS_2}. The batch size is 16,384, and the dimension is 50.}
		\label{fig:ErrorHist_bs2}
	\end{figure}

	For the nonlinear term $\eqref{eq:f_BS_2}$, we also take the region $[a,b]^d=[-0.5,0.5]^d$. Table \ref{tab:bs2_combined} presents the mean and standard deviation of the relative $L^{2}$ error over 8 independent runs, for dimensions 20 and  50, respectively. Both tables demonstrate that the RQMC-based algorithm outperforms the MC-based one. Figure \ref{fig:ErrorHist_bs2} presents the density histogram of pointwise errors between the trained networks and the exact solution. For the MC method, the smallest $L_2$ error is $3.05 \times 10^{-3}$ and the largest $L_2$ error is $1.22 \times 10^{-2}$. For the RQMC method, the smallest $L_2$ error is $1.43 \times 10^{-3}$ and the largest $L_2$ error is $3.88 \times 10^{-3}$. The figure demonstrates that the RQMC-based algorithm achieves smaller errors across the majority of samples.

\begin{table}[H]
    \centering
    \begin{tabular}{|c|c|c|c|c|}
        \toprule
        Dimension & Batch Size & Method & Mean Error & Standard Deviation \\
        \midrule
        \multirow{6}{*}{20} & \multirow{2}{*}{$2^{12}$} & MC & $4.6 \times 10^{-2}$ & $8.8 \times 10^{-3}$ \\
        & & RQMC & $\bm{3.4 \times 10^{-2}}$ & $\bm{1.3 \times 10^{-3}}$ \\
        \cmidrule{2-5}
        & \multirow{2}{*}{$2^{14}$} & MC & $2.0 \times 10^{-2}$ & $7.9 \times 10^{-3}$ \\
        & & RQMC & $\bm{9.3 \times 10^{-3}}$ & $\bm{1.3 \times 10^{-3}}$ \\
        \cmidrule{2-5}
        & \multirow{2}{*}{$2^{16}$} & MC & $9.8 \times 10^{-3}$ & $2.5 \times 10^{-3}$ \\
        & & RQMC & $\bm{5.1 \times 10^{-3}}$ & $\bm{4.0 \times 10^{-4}}$ \\
        \midrule
        \multirow{6}{*}{50} & \multirow{2}{*}{$2^{12}$} & MC & $7.2 \times 10^{-2}$ & $2.1 \times 10^{-2}$ \\
        & & RQMC & $\bm{4.0 \times 10^{-2}}$ & $\bm{2.4 \times 10^{-3}}$ \\
        \cmidrule{2-5}
        & \multirow{2}{*}{$2^{14}$} & MC & $3.1 \times 10^{-2}$ & $7.7 \times 10^{-3}$ \\
        & & RQMC & $\bm{1.1 \times 10^{-2}}$ & $\bm{2.3 \times 10^{-3}}$ \\
        \cmidrule{2-5}
        & \multirow{2}{*}{$2^{16}$} & MC & $1.4 \times 10^{-2}$ & $4.2 \times 10^{-3}$ \\
        & & RQMC & $\bm{6.2 \times 10^{-3}}$ & $\bm{2.3 \times 10^{-3}}$ \\
        \bottomrule
    \end{tabular}
    \caption{The mean and standard deviation of the relative $ L_2$ errors over 8 independent runs for solving the nonlinear Black-Scholes equation with the nonlinear term~\eqref{eq:f_BS_1}.}
    \label{tab:bs1_combined}
\end{table}

        \begin{figure}[H]
		\centering
		\subfloat[Best trained neural networks]{
			\includegraphics[width=0.45\linewidth]{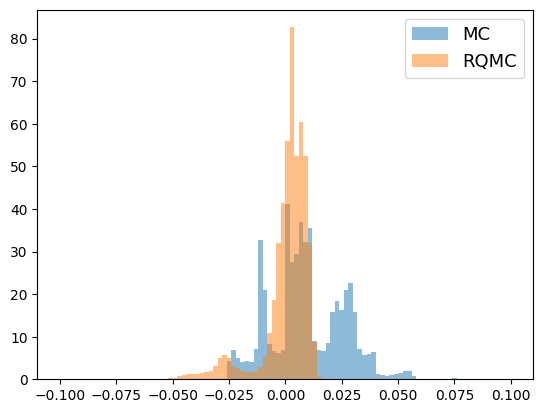}}
		\hspace{2em}
		\subfloat[Worst trained neural networks]{
			\includegraphics[width=0.45\linewidth]{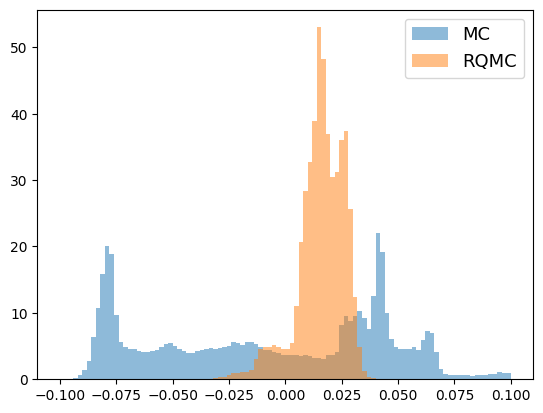}}
		\caption{The density histograms of the pointwise absolute errors for the best and worst trained neural networks in solving the  nonlinear Black-Scholes equation with the nonlinear term \eqref{eq:f_BS_1}. The batch size is 16,384, and the dimension is 50.}
		\label{fig:ErrorHist_bs1}
	\end{figure}

 \begin{table}[H]
    \centering
    \begin{tabular}{|c|c|c|c|c|}
        \toprule
        Dimension & Batch Size & Method & Mean Error & Standard Deviation \\
        \midrule
        \multirow{6}{*}{20} & \multirow{2}{*}{$2^{12}$} & MC & $1.5 \times 10^{-2}$ & $5.7 \times 10^{-3}$ \\
        & & RQMC & $\bm{7.5 \times 10^{-3}}$ & $\bm{1.2 \times 10^{-3}}$ \\
        \cmidrule{2-5}
        & \multirow{2}{*}{$2^{14}$} & MC & $5.7 \times 10^{-3}$ & $4.0 \times 10^{-3}$ \\
        & & RQMC & $\bm{2.4 \times 10^{-3}}$ & $\bm{2.5 \times 10^{-4}}$ \\
        \cmidrule{2-5}
        & \multirow{2}{*}{$2^{16}$} & MC & $3.6 \times 10^{-3}$ & $1.1 \times 10^{-3}$ \\
        & & RQMC & $\bm{2.2 \times 10^{-3}}$ & $\bm{1.4 \times 10^{-4}}$ \\
        \midrule
        \multirow{6}{*}{50} & \multirow{2}{*}{$2^{12}$} & MC & $1.8 \times 10^{-2}$ & $8.2 \times 10^{-3}$ \\
        & & RQMC & $\bm{8.1 \times 10^{-3}}$ & $\bm{2.4 \times 10^{-3}}$ \\
        \cmidrule{2-5}
        & \multirow{2}{*}{$2^{14}$} & MC & $7.8 \times 10^{-3}$ & $3.4 \times 10^{-3}$ \\
        & & RQMC & $\bm{2.4 \times 10^{-3}}$ & $\bm{8.3 \times 10^{-4}}$ \\
        \cmidrule{2-5}
        & \multirow{2}{*}{$2^{16}$} & MC & $3.6 \times 10^{-3}$ & $1.2 \times 10^{-3}$ \\
        & & RQMC & $\bm{1.9 \times 10^{-3}}$ & $\bm{3.3 \times 10^{-4}}$ \\
        \bottomrule
    \end{tabular}
    \caption{The mean and standard deviation of the relative $L^{2}$ errors over 8 independent runs for solving the nonlinear Black-Scholes equation with the nonlinear term~\eqref{eq:f_BS_2}. }
    \label{tab:bs2_combined}
\end{table}

	\section{Conclusion}\label{sec:conclusion}
            Nonlinear PDEs~\eqref{eq:PDE} represent a general form of the linear Kolmogorov PDEs~\cite{beck2021solving,berner2020numerically,de2022error,glau2022deep,richter2022robust}, related to BSDEs by nonlinear Feymann-Kac formula~\cite{pardoux90}, which find extensive applications in  financial mathematics~\cite{karoui97}. Therefore, efficiently solving nonlinear PDEs constitutes a significant challenge.
            We investigated the convergence orders of the generalization error within the DBDP~\cite{hure2020} method under both the MC and QMC approaches. Our theoretical findings indicate that the generalization error under the QMC method exhibits a higher convergence order compared to the MC method. Theoretically, this suggests that the use of the QMC method for sampling during training yields superior outcomes. We validated these theoretical insights through several numerical experiments. The results obtained using the QMC method have smaller errors and improved stability. It is an intriguing prospect to explore the trade-off among scheme error, approximation error, and generalization error to achieve an optimal total error with reduced computational costs. This is left for future work.

\appendix

\section{Proof of Theorem~\ref{thm:total_error}}\label{appendix:total_error} 
	\begin{proof}
		In Theorem 4.1 of Hure et al.~\cite{hure2020}, they propose a mid-term for the proof 
		\begin{equation}\label{eq:def_cV}
			\left\{\begin{aligned}
				\cV_{i}(X_{t_i}) &= \E_i\left[ \cU_{i+1}(X_{t_{i+1}})\right] + f(t_i,X_{t_i},\cV_{i}(X_{t_i}),\bhZ_i(X_{t_i}))\Delta t_i,\\
				\bhZ_i(X_{t_i}) & = \frac{1}{\Delta t_i}\E_i\left[ \cU_{i+1}(X_{t_{i+1}})\Delta W_{t_i}\right],      
			\end{aligned}\right.
		\end{equation}
		where $\E_i\left[ \cdot\right]: = \E\left[ \cdot|X_{t_{i}}\right]$.
		By the martingale representation theorem~\cite{karatzas2014}, there is a $ \hZ_s$ such that 
		\begin{equation}\label{eq:U_i+1}
			\cU_{i+1}(X_{t_{i+1}}) = \cV_i(X_{t_i})-f(t_i,X_{t_i},\cV_i(X_{t_i}),\bhZ_i(X_{t_i}))\Delta t_i + \int_{t_i}^{t_{i+1}} \hZ_s \dd W_s.
		\end{equation}
		and by It\^o isometry, we have 
		\begin{equation}
			\bhZ_i(X_{t_i}) = \frac{1}{\Delta t_i} \E_{i}\left[ \int_{t_i}^{t_{i+1}}\hZ_s \dd s\right].
		\end{equation}
		Hure et al.~\cite{hure2020} obtain the following result
		\begin{equation}\label{eq:error_initial}
			\begin{aligned}
				\varepsilon\left[ \left( \cU,\cZ\right),\left( Y,Z\right)\right] \le C\bigg( &\E\left| g(\X_T) - g(X_N)\right|^2 + |\pi| + \varepsilon^{Z}(\pi) \\
				&+ \sum_{i =0}^{N-1}\Big(N\E\left| \cV_i(X_{t_i}) - \cU_i(X_{t_i})\right|^2 + \E\left| \bhZ_i(X_{t_i}) - \cZ_i(X_{t_i})\right|^2\Big)\bigg),
			\end{aligned}
		\end{equation}
		where $$\varepsilon^{Z}(\pi) = \E\left[ \sum_{i =0}^{N-1}\int_{t_i}^{t_{i+1}}\left|Z_t - \overline{Z}_{t_i} \right|^2\dd t\right],\quad  \overline{Z}_{t_i} = \frac{1}{\Delta t_i} \E_i\left[\int_{t_i}^{t_{i+1}}Z_t\dd t\right].$$
		Note that by substituting~\eqref{eq:U_i+1} into $L_i(\cU_{i+1},\cU_i,\cZ_i)$, we can decompose it as follows
		\begin{align}\label{eq:decomp_L}
			L_i(\cU_{i+1},\cU_i,\cZ_i) = \tilde L_i(\cU_i,\cZ_i) + \E\left[ \int_{t_i}^{t_{i+1}}\left| \widehat Z_t - \bhZ_i(X_{t_i})\right|^2\dd t\right],
		\end{align}
		where 
		\begin{equation}
			\begin{aligned}
				&\tilde L_i(\cU_i,\cZ_i)= \Delta t_i \E\left| \bhZ_i(X_{t_i}) - \cZ_i(X_{t_i})\right|^2\\
				+ &\E\left| \cV_i(X_{t_i}) - \cU_i(X_{t_i}) + \left(f(t_i,X_{t_i},\cU_i(X_{t_i}),\cZ_i(X_{t_i})) - f(t_i,X_{t_i},\cV_i(X_{t_i}),\bhZ_i(X_{t_i}))\right)\Delta t_i\right|^2.
			\end{aligned}
		\end{equation}
		Note that the second term of the right hand side of~\eqref{eq:decomp_L} is independent of $(\cU_i,\cZ_i)$, thus it is the minimize of $L_i(\cU,\cZ)$. By~\eqref{eq:def_cV}, $\cV$ and $\bhZ$ are functions of $X_{t_i}$. Therefore, the minimizer $(\cU^{*},\cZ^{*})$ of~\eqref{eq:decomp_L} satisfies
		\begin{equation}\label{eq:uz_optimal}
			\cU_i^*(X_{t_i}) = \cV_{t_i}(X_{t_i}),\quad \cZ_i^*(X_{t_i}) = \bhZ_i(X_{t_i}),
		\end{equation}
		Moreover,
		\begin{equation}
			\tilde L_i(\cU_i^*,\cZ_i^*) =  0,\quad L_i(\cU_i^*,\cZ_i^*) = \E\left[ \int_{t_i}^{t_{i+1}}\left| \hZ_t - \bhZ_i(X_{t_i})\right|^2 \dd t\right].
		\end{equation}
		By the lipschitz condition of $f$, there is a constant $C$, such that
		\begin{equation}\label{eq:estimate_tilde_L}
			\left\{\begin{aligned}
				\tilde L_i\left( \cU_i,\cZ_i\right) &\le \left( 1 + C\Delta t_i\right)\E\left| \cV_i(X_{t_i}) - \cU_i(X_{t_i})\right|^2 + C\Delta t_i\E\left| \bhZ_i(X_{t_i}) - \cZ_i(X_{t_i})\right|^2,\\
				\tilde L_i\left( \cU_i,\cZ_i\right) &\ge \left( 1 - C\Delta t_i\right)\E\left| \cV_i(X_{t_i}) - \cU_i(X_{t_i})\right|^2 + C\frac{\Delta t_i}{2}\E\left| \bhZ_i(X_{t_i}) - \cZ_i(X_{t_i})\right|^2.
			\end{aligned}\right.
		\end{equation}
		It follows from~\eqref{eq:error_initial},~\eqref{eq:decomp_L} and~\eqref{eq:estimate_tilde_L} that
		\begin{equation}
			\begin{aligned}
				\varepsilon\left[ \left( \cU,\cZ\right),\left( Y,Z\right)\right] \le C\bigg( &\E\left| g(\X_T) - g(X_N)\right|^2 + |\pi| + \varepsilon^{Z}(\pi) \\
				&+ N\sum_{i =0}^{N-1}\left( L_i\left( \cU_{i+1},\cU_i,\cZ_i\right) - L_i\left( \cU_{i+1},\cU_i^*,\cZ_i^*\right)\right)\bigg),
			\end{aligned}
		\end{equation}
		where $\left( \cU_i^*,\cZ_i^*\right) = \arg\min_{(\cU,\cZ)} L_i\left( \cU_{i+1},\cU,\cZ\right)$, which has the representation~\eqref{eq:uz_optimal}. Noting that $g$ has the lipschitz condition, by Zhang~\cite{zhang04}, we have
		\begin{equation}\notag
			\varepsilon^Z(\pi) = O(|\pi|).
		\end{equation}
		Thus under the conditions of Theorem~\ref{thm:total_error},~\eqref{eq:total_error_final} holds.
	\end{proof}

\end{document}